\documentclass[11pt,onecoulme]{article}
\usepackage{amsfonts}
\usepackage{mathrsfs}
\usepackage{amssymb}
\usepackage{color, amsmath,amssymb, amsfonts, amstext,amsthm, latexsym}

\usepackage{amssymb, epsfig, amssymb, latexsym}
\usepackage{amsmath}
\usepackage{graphicx}
\usepackage{longtable}

\numberwithin{equation}{section}

\allowdisplaybreaks

\newtheorem{definition}{Definition}[section]
\newtheorem{theorem}[definition]{Theorem}
\newtheorem{lemma}[definition]{Lemma}
\newtheorem{remark}[definition]{Remark}

\newtheorem{hyp}[definition]{Hypothesis}
\title{ A Class of Delay Optimal Control Problems  and Viscosity
           Solutions to
           Associated
               Hamilton-Jacobi-Bellman Equations \thanks{This work was partially supported by  the National Natural Science Foundation of China  (Grant No. 11401474),  Shaanxi Natural Science Foundation
               (Grant No. 2014JQ1035)
              and the Fundamental Research Funds for the Central Universities (Grant No. 2452015087).}}

\author{Jianjun Zhou 
 \\
College of Science,
             Northwest A\&F University,\\ Yangling 712100, Shaanxi, P. R.
             China\\
      \emph{E-mail:zhoujj198310@163.com}}

\begin{document}

\maketitle

\pagestyle{plain}

\begin{abstract}
 In this article,  a class of  optimal control problems of differential equations with
         delays are investigated  for which the associated Hamilton-Jacobi-Bellman (HJB)
          equations are nonlinear partial differential equations
          with delays.
         %, where the state equations depend
%       %  not only on the state of a certain point in the past, but also on the state of a certain period in the past.
%        The
%         problems are considered in a  bounded, right continuous functions space ${\cal{D}}$ and  the associated  Hamilton-Jacobi-Bellman (HJB)
%          equations are studied.
       This type of HJB equation has not been previously  studied  and is difficult to solve because the state equations do not possess
       smoothing properties.  We introduce a slightly different notion of viscosity
            solutions and
             identify the value function of the optimal  control problems as
             a
             unique viscosity solution to the associated HJB equations. 
\medskip

 {\bf Key Words:}  Hamilton-Jacobi-Bellman equations; Viscosity
                 solutions; Optimal control;
                  Differential equations with delays; Existence
        and uniqueness
\end{abstract}

{\bf 2000 AMS Subject Classification:} 34K35; 49L20; 49L25.

\section{Introduction}
%%%%%%%%%%%%%%%%%%%%%%%%%%%%%%%%%%%%%%%%%%%%%%%%%%%%%%%%%%%%%%%%%

\par
In this paper, we consider the following   controlled   differential
                 equations with delays:
\begin{equation}
\begin{cases}
dX^u(s)=
            F{(}s,X^u(s),(a,X^u_s)_{H},u(s){)}ds+b(s)X^u(s-\tau)ds,   s\in [t,T],\\
 ~~~~~X^u_t=x\in
        {\cal{D}},
\end{cases}
\end{equation}
%$$
%            \cases{dX^u(s)=
%            F{(}s,X^u(s),(a,X^u_s)_{H},u(s){)}ds+b(s)X^u(s-\tau)ds,   s\in [t,T],\cr
%             ~~~~~X^u_t=x\in
%        {\cal{D}},
%           }
%           \eqno(1.1)
%$$
        where
$$
        X^u_s(\theta)=X^u(s+\theta),\ \theta\in[-\tau,0],\ (a,X^u_s)_{H}=\int^{0}_{-\tau}(a(\theta),X^u_s(\theta))_{R^d}d\theta.
$$
%
%$$
%            \cases{dX^u(s)=
%            F(s,X^u(s),X^u(s-\tau),X^u_s,u(s))ds,   s\in [t,T],\cr
%            ~~~~~X^u_t=x,}
%           \eqno(1.1)
%$$
%        where
%$$
%        X^u_s(\theta)=X^u(s+\theta),\ \theta\in[-\tau,0], \ \
%          x\in {\cal{D}}.
%$$
                In the equations above, ${\cal{D}}$ denotes the space of bounded, right continuous, $R^d$-valued
         functions on $[-\tau,0]$, and $F:[0,T]\times R^d\times R\times U\rightarrow R^d$ is a given
         map,
         where  $U$ is a  metric space in which the control $u(\cdot)$ takes values. For any initial state $x\in
         {\cal{D}}$ and control $u(\cdot)\in {\cal{U}}[t,T]:=\{u:[t,T]\rightarrow U|\ u(\cdot) \mbox{ is
         measurable}\}$, the corresponding trajectory $X(\cdot)$ is
         a solution to (1.1).
               % The unknown $X^u(s)$, representing the state of the system, is a measurable function from $[t,T]$ in $R^d$;
%                the control process $u$ takes values in a  metric space $U$.
$a$ and $b$ %$a\in W^{1,2}([-\tau,0];R^d),\ b\in W^{1,2}([0,T];R^{d\times d})$
                are two given
                functions that  satisfy suitable  smoothness
                properties,
                 and the coefficient $F$ is assumed to satisfy a Lipschitz condition with
        respect to the appropriate norm.  Thus, the solution to (1.1) is uniquely determined by the initial state  and the control.
\par
               The control problem consists of minimizing a cost functional of
               the following
            form:
\begin{equation}
   J(t,x,u)=\int_{t}^{T}q(\sigma,X^u(\sigma),u(\sigma))d\sigma+\phi(X^u(T)), 
\end{equation}
            over all of the  controls $u(\cdot)\in {\cal{U}}[t,T]$.
            Here, $q$ and $\phi$ are functions on $[0,T]\times R^d\times U$ and $R^d$, respectively.
                           We define the value function of the  optimal
                  control problem as follows:
\begin{equation}
                  V(t,x):=\inf_{u\in{\mathcal
                  {U}}[t,T]}J(t,x,u), \ \ \
                   t\in [0,T],\ x\in {\cal{D}}.
\end{equation}
                 We  assume that $q$ and $\phi$ satisfy suitable
                     conditions and consider  the following Hamilton-Jacobi-Bellman
                     (HJB) equations:
\begin{equation}
\begin{cases}
\frac{\partial}{\partial t} V(t,x)+{\mathcal
                        {S}}(V)(t,x)+ H(t,x,\nabla_x V(t,x))= 0, \ \  t\in
                               [0,T],\ \ x\in {\mathcal{D}},\\
 V(T,x)=\phi(x(0)),
\end{cases}
\end{equation}
%$$
%          \cases{\frac{\partial}{\partial t} V(t,x)+{\mathcal
%                        {S}}(V)(t,x)+ H(t,x,\nabla_x V(t,x))= 0, \ \  t\in
%                               [0,T],\ \ x\in {\mathcal{D}}, \cr
%          V(T,x)=\phi(x(0)),\cr}\eqno(1.4)
%$$
                where
\begin{eqnarray*}
                                H(t,x,p)=\inf_{u\in{
                                         {U}}}[\langle
                        [F(t,x(0),(a,x)_{H},u)+b(t)x(-\tau)]1_0(t),p\rangle+q(t,x(0),u)], (t,x,p)\in [0,T]\times {\cal{D}}\times
                        {\cal{D}}^*.
\end{eqnarray*}
                       Here $1_0$ denotes the character function of {\{0\}}. The definition of a weak infinitesimal generator
                         $\cal{S}$ will be given in section 3.
\par
                         The primary objective of this article is to
                         develop the notion of a viscosity solution
                         to
                         the HJB equations given by (1.4). We then show the value function
                         $V$  defined in  (1.3) is a unique viscosity solution to the  HJB
                         equations given in  (1.4).
\par
               The type  of problem above arises
               in  many different fields of application, including  engineering, economics and
               biology. These problems  typically disrupt  the
               optimum operation of a system in the form of a time lag in the response to a given
              input.  References  \cite{[1]}, \cite{bam}, \cite{bou}, and  \cite{fab} present
              models with delays  in economics; references
              \cite{car} and  \cite{fab1} present deterministic advertising models with
              delay effects; references  \cite{bou1} and \cite{bou2} present  population
              models.
\par
              These optimal control problems for differential  equations with
              delays have been thoroughly investigated in recent years
              (see \cite{bou}, \cite{bou1}, \cite{bou2},
              \cite{car}, \cite{car1}, \cite{fab}, \cite{fab1}, and  \cite{fed}). However, to the best of our knowledge, none of these results
              are directly applicable to our case.  In reference \cite{car1} and \cite{fed},  the optimal control problem was embedded in a Hilbert
              space,
                and  the viscosity solutions for the associated HJB
              equations were investigated. These results  do not hold
              when
              $b(\cdot)\neq0$ in the state equation. In  \cite{bou1}, \cite{bou2}, \cite{fab}, and
              \cite{fab1},
              the term
              $b(\cdot)X(\cdot-\tau)$ in the  state equation is considered,  the results obtained in these references only apply when   the
              state equation is a linear differential equation with
              delays.  Optimal control
              problems of a state equation with memory were investigated in \cite{car}; in these problems,  however, the control must satisfy a linear condition,
               which  is not fulfilled if $F$ is a genuinely
              nonlinear function.
\par
              %The optimal control problem (1.1) and (1.2) has not been studied
%              in detail (see, e.g. [4, 13]).
          It is well known that the optimal control problem given by (1.1) and (1.2)  can be reformulated as an optimal control problem of the evolution equation
                in a Hilbert space  (see, e.g.
                \cite{car1} and \cite{fed}).
                 %Federico,  Goldys, and Gozzi [13] embedded the optimal control problem in a Hilbert space
%                and  studied the regularity of viscosity solutions for the associated HJB
%              equations. However, the results in [13] are not applicable to our situation. In fact,
                However, in this case,  the initial value must have the following form:
                 $X^u_t=x\in H$ and $X^u(t)=x^0\in R^d$. This form ensures that  the value
                function is not a viscosity solution of the associated
                HJB equations because of the   $b(\cdot)X^u(\cdot-\tau)$ term  in the state
                equations. To the aforementioned  challenges, we  study the
                associated HJB equations in an infinite dimensional
                space ${\cal{D}}$.
\par
             Crandall-Lions
           \cite{cra1}introduced the notion of viscosity solutions to the HJB equations in the early
           1980's and
           showed that  the dynamic programming method could  be applied to optimal control problems. Since then, many papers have been published on the
           development of
            the theory of viscosity solutions (see, e.g., \cite{bar}, \cite{cra}, \cite{cra2}, \cite{ish}, \cite{ish1}, and \cite{yon}).
             References for dealing with equation in  an
                  infinite dimensional Hilbert space include  \cite{cra3}, \cite{cra4}, \cite{cra5}, \cite{cra6}, \cite{cra7}, and \cite{li}.
              In references \cite{cra3}, \cite{cra4}, \cite{cra5}, \cite{cra6}, and  \cite{cra7},
               Crandall  and Lions systematically introduced  the basic theories for  viscosity
              solutions.  Then, Zhou and Yong \cite{li} proved  the existence and uniqueness of
              a viscosity solution to general unbounded first-order HJB equations in infinite dimensional Hilbert spaces.
%\par
%              The notion of viscosity solution has been
%              successfully applied to HJB equaition and subsequently
%              to optimal control problems.
\par
              To the  best
                of
             our knowledge, the associated HJB equations (1.4) have  not
             been previously
                studied. The primary difficulty in solving these equations is caused by the infinite-dimensionality of the space of variables
                and thus the non-compactness of the space. Hence,our
                problem
                does not fall into the framework used in references   \cite{cra3}, \cite{cra4}, \cite{cra5}, \cite{cra6}, \cite{cra7}, and \cite{li}.
                Thus,
                the standard proofs of the comparison theorem rely heavily on  compactness
                arguments and are not applicable to
                 our case.
\par
             To overcome this difficulty,  we first  prove a left
              maximization principle for the   space $[0,T]\times {\cal{D}}\times[0,T]\times {\cal{D}}$ ( see Lemma
              4.1), i.e., variables exist that maximize functions defined. The proof of the comparison
              theorem involves  maximizing the
              auxiliary function. The underlying principle  is to use the left maximization principle to find a
              variable that maximizes the auxiliary
              function.
\par
             We next   introduce a  slightly different  notion of  a viscosity solution to the HJB equations given in  (1.4).  We use the left
              maximization principle to prove  the uniqueness of a viscosity
                 solution that corresponds to  our new definition og s viscosity solution. At the same time,
                 we show that
                 the value function is  a viscosity solution to
                 the HJB equations.
\par
                    Our results rely heavily on
                      the construction of state equations. We hope to overcome this
                     serious limitation of our approach in  future work. However,
                     our method is suitable for a large class of
                     optimal control problems for differential
                     equations with delays. 
\par
                 The paper is organized as follows. In the following
              section, we define our notation and review the background for  differential equations with delays are studied.
              In section 3 we prove   the dynamic programming principle (DPP) and Lemma 3.5 which are used in the following sections.
               In section 4, we define viscosity solutions and  show that the value function $V$ defined by (1.3) is a viscosity solution
                to the HJB equations given in  (1.4).
                Finally,  the uniqueness of viscosity solutions to (1.4) is proved in section 5.

\section{Preliminary work}  %\label{RDS}
%%%%%%%%%%%%%%%%%%%%%%%%%%%%%%%%%%%%%%%%%%%%%%%%%%%%%%%%%%%%%%%%%

\par
   Here, we define the notations that are used in this paper.
   %The norm of an element $x$ of a Banach space $F$ will be denoted $|x|_F$ or simply $|x|$, if no confusion is possible.
    We use the symbol $|\cdot|$ to denote the norm in
         a Banach space $F$, the norm symbol is subscripted when necessary. For the vectors $x,y\in R^d$, the scalar product is denoted by $(x,y)_{R^d}$ and the
         Euclidean norm $(x,x)^{\frac{1}{2}}_{R^d}$ is denoted by $|x|$. For $T>0$ and $0\leq t<T$, let $C([t,T],R^d)$
        denote
        the space of continuous functions from $[t,T]$ to $R^d$, which is associated with the usual
         norm $|f|_C=\sup_{{\theta\in [t,T]}}|f(\theta)|$. Let $\tau>0$ be fixed; then    $H$ denotes  %and $\mathbf{H}$ % and $\cal{H}$
         the real, separable Hilbert space $L^2([-\tau,0];R^d)$ %and $R^d\times L^2([-\tau,0];R^d)$, %and $R^d\times R^d\times L^2([-\tau,0];R^d)$,
          for scalar product
          $(\cdot,\cdot)_H$. 
          Let
         ${\cal{D}}$ denote the set of bounded, right continuous, $R^d$-valued
         functions on $[-\tau,0]$. We define a norm on $\cal{D}$  as
         follows:
   $$
                      |\omega|_{\cal{D}}= \sup_{\theta\in
                      [-\tau,0]}|\omega(\theta)|,\ \ \omega \in {\cal{D}}.
   $$
               Then, $({\cal{D}},|\cdot|_{\cal{D}})$ is  Banach space.
               \par
              We define the $|\cdot|_{B}$-norm on $H$ as follows:
$$
                              |x|^2_{B}:=\int^{0}_{-\tau}(Bx)^2(s)ds,
$$
              where
 $$
            ( Bx)(s)=\int^{0}_{s}x(\theta)d\theta,\ \ s\in [-\tau,0].
 $$
\par
         Let $0\leq t\leq \bar{t}\leq T$, $0\leq s\leq \bar{s}\leq T$, and $\omega, \bar{\omega},\nu, \bar{\nu}\in
         {\cal{D}}$ be given. We define  $({t},{\omega})\otimes (\bar{t},\bar{\omega})
                         \in [0,T]\times {\cal{D}}$  and
 $({t},{\omega},s, \nu)\otimes (\bar{t},\bar{\omega},\bar{s},\bar{\nu})
                         \in [0,T]\times {\cal{D}}\times[0,T]\times {\cal{D}}$  by
$$
                ({t},{\omega})\otimes (\bar{t},\bar{\omega}):=
                (\bar{t},\tilde{{\omega}}),\ \ \ \ \
                ({t},{\omega},s, \nu)\otimes
                (\bar{t},\bar{\omega},\bar{s},\bar{\nu}):=
                (\bar{t},\tilde{{\omega}},\bar{s},\tilde{{\nu}}),
$$
where
\begin{equation*}
\tilde{{\omega}}(\theta)=
\begin{cases}
 \bar{\omega}(\theta), \ \ \ \ \  \ \ \  \ \ \ \  \ t-\bar{t}\leq\theta\leq0,\\
 {\omega}( \bar{t}-t+\theta), \ \ \
                                                         -\tau\leq\theta<t-\bar{t},
\end{cases}
\end{equation*}
\begin{equation*}
 \tilde{\nu}(\theta)=
\begin{cases}
  \bar{\nu}(\theta), \ \ \ \ \  \ \ \ \ \ \ \  \ t-\bar{t}\leq\theta\leq0,\\
 \nu( \bar{t}-t+\theta), \ \ \ -\tau\leq\theta<t-\bar{t}.
\end{cases}
\end{equation*}
       %where
%$$
%            \tilde{{\omega}}(\theta)=\cases{ \bar{\omega}(\theta), \ \ \ \ \  \ \ \  \ \ \ \  \ t-\bar{t}\leq\theta\leq0,\cr
%                                                         {\omega}( \bar{t}-t+\theta), \ \ \
%                                                         -\tau\leq\theta<t-\bar{t},\cr}\
%                                                         \ \ \ \
%            \tilde{\nu}(\theta)=\cases{ \bar{\nu}(\theta), \ \ \ \ \  \ \ \ \ \ \ \  \ t-\bar{t}\leq\theta\leq0,\cr
%                                                         \nu( \bar{t}-t+\theta), \ \ \ -\tau\leq\theta<t-\bar{t}.\cr}
%$$
%$$
%         ({t},{\omega}(\theta))\otimes (\bar{t},\bar{\omega}(\theta)):= \cases{ (\bar{t},\bar{\omega}(\theta)), \ \ \ \ \  \ \ \ \ \  \ t-\bar{t}\leq\theta\leq0,\cr
%                                                          (\bar{t},{\omega}( \bar{t}-t+\theta)), \ \ \ -\tau\leq\theta<t-\bar{t}.\cr}
%$$
We denote the boundary of a given open subset $Q\subset R^d$ by
$\partial Q$ and $\bar{Q}=Q\bigcup \partial Q$. Let us define
$$
                         {\cal{D}}_Q:=\{\omega\in {\cal{D}}: \omega(\theta)\in Q,\ \theta\in
                                        [-\tau,0]\}
$$
   and
$$
                         {\cal{D}}_{\bar{Q}}:=\{\omega\in {\cal{D}}: \omega(\theta)\in \bar{Q},\ \theta\in
                         [-\tau,0]\}.
$$
            % We will say that a measurable process $u$
%             with values in a  metric space $U$ is an $admissible \ control$, we denote by $\mathcal {U}$
%              the set of admissible controls.
               Let us consider the controlled state
             equations:
\begin{equation}
\begin{cases}
dX^u(s)=
            F{(}s,X^u(s),(a,X^u_s)_{H},u(s){)}ds+b(s)X^u(s-\tau)ds, \  s\in [t,T],\\
 ~~~~~X^u_t=x\in {\cal{D}},
\end{cases}
\end{equation}
%$$
%            \cases{dX^u(s)=
%            F{(}s,X^u(s),(a,X^u_s)_{H},u(s){)}ds+b(s)X^u(s-\tau)ds, \  s\in [t,T],\cr
%            ~~~~~X^u_t=x\in {\cal{D}},
%           }
%           \eqno(2.1)
%$$
        where
$$
        X^u_s\in {\cal{D}},\ \  X^u_s(\theta)=X^u(s+\theta),\ \theta\in[-\tau,0].%,\ (a,X^u_s)_{R^d}=\int^{0}_{-\tau}a(\theta)X^u_s(\theta)d\theta,\
$$
   Here, the control $u(\cdot)$ belongs to
$$
                             {\cal{U}}[t,T]:=\{u(\cdot):[t,T]\rightarrow U|\ u(\cdot) \ \mbox{is
                             measurable}\},
$$
        and where  $U$ is a metric space.
              We make the
           following assumptions.
\begin{hyp}
\begin{description}
        \item{(i)}
                 The mapping $F$: $[0,T]\times  R^d\times  R\times U\rightarrow R^d$
                   is measurable and  a constant
                           $L>0$ exists  such that, for every $t,s\in [0,T], x,y\in R^d\times R,   u\in U$,
\begin{eqnarray*}
                |F(t,x,u)|\leq L(1 + |x|)\ \ \ \mbox{and} \ \ \
                  |F(t,x,u)-F(s,y,u)|\leq
                 L(|s-t|+|x-y|).
\end{eqnarray*}
\item{(ii)}
    $a(\cdot)\in W^{1,2}([-\tau,0];R^d)$ with  $a(-\tau)=0$ and $b(\cdot)\in
    W^{1,2}([0,T];R^{d\times d})$, and  a constant $L>0$ exists  such
    that, for every $t,s\in [0,T]$,
    $$
                         |b(s)-b(t)|\leq L|s-t|.
    $$
                   %and for every $(t,x,y,z,u)\in [0,T]\times R^d\times H\times U$,
%$|F(t,x_1,y_1,z_1)-F(s,x_2,y_2,z_2)|\leq
%                 L(|s-t|+|x_1-x_2|+|y_1-y_2|+|z_1-z_2|_B),\ t, s\in [0,T], \ \ x_1,x_2,y_1,y_2\in R^d, z_1,z_2 \in H.
%$
%\\
% $ |F(t^1,x^1,u)-F(t^2,x^2,u)|\leq
%                 L|x^1-x^2|+\gamma(|t^1-t^2|,(|x^1|+|z^1|)),
%$\\
%$
%               |F(t,x,u)|\leq L(1 + |x|),   \ t^1,t^2\in [0,T], \ u
%                 \in U,\ \ x, x^1,x^2\in {\cal{D}}.
%$
\end{description}
\end{hyp}
\par
              A function  $X^u:[t,T]\rightarrow R^d$ is a  solution to equation $(2.1)$ if  the function satisfies the following condition:
$$
             X^u(s)=x(0)+\int_{t}^{s}F{(}\sigma,X^u(\sigma),(a,X^u_\sigma)_{H},u(\sigma){)}d\sigma
                    +\int_{t}^{s}b(\sigma)X^u(\sigma-\tau)d\sigma,\ \
             s\in [t,T], \eqno(2.2)
$$
       where $X^u_t=x\in {\cal{D}}, X^u(s)=x(s-t), \ t-\tau\leq s<t$.
              To emphasize the  dependence of the solution on the initial data, we denote the solution by
                $X^u(s,t,x)$.
\par
\begin{theorem} \ \ Let us assume that Hypothesis 2.1  holds. Then,
a unique function
             $X\in C([t,T];R^d)$ exists that is a solution to
             $(2.1)$. Moreover,
$$
                  \sup_{s\in [t,T]}|X^{u}(s,t,x)|\leq C_1\bigg{(}1+|x(0)|+\sup_{l\in[-\tau,0]}
                          \bigg{|}\int^{l}_{-\tau}x(\theta)d\theta\bigg{|}+|x|_{B}\bigg{)}\leq C_2{(}1+|x(0)|+|x|_H{)},\eqno(2.3)
$$
              where the constants $C_1$ and $C_2$ depend only on $L$, $T$, $\tau$ $a(\cdot)$ and $b(\cdot)$.
\end{theorem}
\par
{\bf  Proof}. \ \
            For every initial value $x\in {\cal{D}}$, we define the
                       mapping $\Phi$ from $ C([t,T];R^d)$
                       to itself as
\begin{eqnarray*}
                       \Phi(X^u)(s)=x(0)+\int_{t}^{s}
                        F{(}\sigma,X^{u}(\sigma),(a,X^u_\sigma)_{H},u(\sigma){)}d\sigma
             +\int_{t}^{s}b(\sigma)X^{u}(\sigma-\tau)d\sigma,\ s\in
                       [t,T],
\end{eqnarray*}
                 where   $X^u(s)=x(s-t)$ if  $s<t$. We first show
                 that $\Phi(X^u)$ is continuous with respect to the
                 time $s$. To this end, for every $t\leq s_1\leq s_2\leq
                 T$, there is a  constant $C>0$ that satisfies the following
                 condition:
\begin{eqnarray*}
                  |\Phi(X^u)(s_1)-\Phi(X^u)(s_2)|&\leq& L\int_{s_1}^{s_2}
                        (1+|X^{u}(\sigma)|+|(a,X^u_\sigma)_{H}|{)}d\sigma
             +\int_{s_1}^{s_2}|b(\sigma)||X^{u}(\sigma-\tau)|d\sigma\\
             &\leq&C(1+\sup_{\sigma\in[t,T]}|X^{u}(\sigma)|+|x|_{\cal{D}})|{s_2}-{s_1}|.
\end{eqnarray*}
                     We next show that it is a contraction, under an equivalent norm.  We define the norm
                       ${\parallel X^u\parallel}={\sup}_{s\in[t,T]}e^{-\beta
                       s}|X^u(s)|$, where $\beta>0$ will be chosen later.
                       This norm is equivalent to the original norm on the space $ C([t,T];R^d)$.
                       Then, the definition of the mapping yields
\begin{eqnarray*}
                        &&\parallel \Phi(X^u)\parallel=\sup_{s\in [t,T]}|e^{-\beta s}\Phi(X^u)(s)|\\
                        &\leq&|x(0)|+\sup_{s\in [t,T]}e^{-\beta
                        s}\bigg{[}
                       \int_{t}^{s}
                        {|}F{(}\sigma,X^{u}(\sigma),(a,X^{u}_\sigma)_H,u(\sigma){)} {|}d\sigma
             + \bigg{|}\int_{t}^{s}b(\sigma)X^{u}(\sigma-\tau)d\sigma \bigg{|}\bigg{]}\\
                       %&\leq& |x(0)|+\sup_{s\in [t,T]}e^{-\beta s}\bigg{[}\int^{s}_{t}
%                       L(1+|X^u(\sigma)|+|a|_{W^{1,2}}|X^u_\sigma|_B)d\sigma+|b(s)|\int^{s}_{t}|X^{u}(\sigma)|d\sigma\\
%                        &&       +|b(s)|\sup_{l\in[-\tau,0]}\bigg{|}\int^{l}_{-\tau}x(\sigma)d\sigma\bigg{|}
%                        +|b|_{W^{1,2}}(T-t)^{\frac{1}{2}}\bigg{(}\sup_{l\in[-\tau,0]}\bigg{|}\int^{l}_{-\tau}x(\theta)d\theta\bigg{|}
%                        +\int^{s}_{t}|X^{u}(\sigma)|d\sigma\bigg{)}\bigg{]}\\
                      % &\leq& |x(0)|+(T-t)\bigg{[}L+[\sup_{s\in [0,T]}|b(s)|+|b|_{W^{1,2}}(T-t)^{\frac{1}{2}}]\sup_{l\in[-\tau,0]}
%                          \bigg{|}\int^{l}_{-\tau}x(\theta)d\theta\bigg{|}
%                       +L|a|_{W^{1,2}}|x|_B\bigg{]}\\
%                       &&+(L+\sup_{s\in [0,T]}|b(s)|+|b|_{W^{1,2}}(T-t)^{\frac{1}{2}}+2L\tau |a|_{W^{1,2}})\sup_{s\in [t,T]}\bigg{(}\int^{s}_{t}
%                        e^{-\beta (s-\sigma)}\sup_{r\in[t,\sigma]}e^{-\beta r }|X^u(r)|)d\sigma\bigg{)}\\
                        &\leq& |x(0)|+(T-t)\bigg{[}L+[\sup_{s\in [0,T]}|b(s)|+|b|_{W^{1,2}}(T-t)^{\frac{1}{2}}]\sup_{l\in[-\tau,0]}
                          \bigg{|}\int^{l}_{-\tau}x(\theta)d\theta\bigg{|}
                       +L|a|_{W^{1,2}}|x|_B\bigg{]}\\
                       &&+\frac{1}{\beta}(L+\sup_{s\in [0,T]}|b(s)|+|b|_{W^{1,2}}(T-t)^{\frac{1}{2}}+2L\tau |a|_{W^{1,2}})\parallel X^u\parallel.
                       \ \ \ \ \ \ \ \ \ \ \ \ \ \ \ \ \ \ \ \ \ \ \ \
                         \ \ \ \ \  \ \ \ \ \ \  (2.4)
\end{eqnarray*}
             This result shows that $\Phi$ is a well-defined mapping on $C([t,T];R^d)$. If $X^u,\ X^u_1$ are functions belonging
                 to this space, similar sequences of inequalities  show that
$$
                        \parallel \Phi(X^u)-\Phi(X^u_1)\parallel
                        \leq \frac{1}{\beta}(L+\sup_{s\in [0,T]}|b(s)|+|b|_{W^{1,2}}(T-t)^{\frac{1}{2}}+2L\tau |a|_{W^{1,2}})\parallel X^u-X^u_1\parallel.
                        \eqno(2.5)
$$
                      Therefore, for a sufficiently large $\beta$, the mapping $\Phi$
                     is a contraction.  In addition,  (2.4) can be used to obtain
                     (2.3). This result completes the proof. \ \  $\Box$
                         
\begin{remark}
\begin{description}
\rm{
        \item{(i)}
       The theorem above show that the solution $X^u(\cdot)$ to equation
        (2.1) is continuous with respect to the time $s\in [t,T]$  even if  the initial value $x$ belongs to ${\cal{D}}$.
       \item{(ii)}  Theorem 2.2 also holds true when the initial state $X^u_t=x\in {\cal{D}}$ is
        replaced
        by $X^u_t=x\in H$ and $X^u(t)=x^0\in R^d$.}
\end{description}
\end{remark}
\par
               Let us now consider  some continuities of the solution $X^u(\cdot)$ to equation
        (2.1), these properties will be used in the proof of Theorem 3.2.
\begin{theorem}\ \ Let us assume that Hypothesis 2.1  holds. Then,  constants $C_3, C_4>0$ exist that depend only on $L$, $T$, $\tau$
$a(\cdot)$ and $b(\cdot)$,
   such that, for every  $t,t_1,t_2\in [0,T]$, and  $x_1,x_2\in{\cal{D}}$,
\begin{eqnarray*}
            &&\sup_{u\in {\cal{U}}[t_1\wedge t_2,T]}\sup_{s\in[t_1\vee t_2,T]}|X^{u}(s,t_1,x_1)-X^{u}(s,t_2,x_2)|
              \leq C_3(1+|x_1(0)|+|x_2(0)|_H+|x_1|+|x_2|_H)\\
              &&~~~~~~~~~~~~~~~~~~~\times\bigg{(}|x_1(0)-x_2(0)|+|t_2-t_1|^{\frac{1}{2}}
                  +\sup_{l\in[-\tau,0]}\bigg{|}\int^{0}_{l}x_{1}(\theta)-x_{2}(\theta)d\theta\bigg{|}\bigg{)}, \  \ \ \ \ \ \ \
                   \  \ \ \  (2.6)
\end{eqnarray*}
$$
             \sup_{u\in {\cal{U}}[t,T]}\sup_{s\in[t,T]}|X^{u}(s,t,x_1)-X^{u}(s,t,x_2)|\leq
             C_4(|x_1(0)-x_2(0)|+|x_1-x_2|_H).\eqno(2.7)
$$
\end{theorem}
   \par
{\bf  Proof}. \ \
                    For any $t_1,t_2\in [0,T]$ and $x_1,x_2\in{\cal{D}}$, we  assume
                that $t_1\leq t_2<t_1+\tau$. Let $X^{u,i}(s)$ denote $X^u(s,t_i,x_i)$ for $s\in [t_i,T]$, where $i=1,2$. Thus, we
                obtain the following results:
\begin{eqnarray*}
                  && \sup_{s\in[t_1\vee t_2,l]}|X^{u,1}(s)-X^{u,2}(s)|\\
                %&\leq&
%                |x_1(0)-x_2(0)|+\int_{t_1}^{t_2}{|}F(\sigma,X^{u,1}(\sigma),
%                       (a,X^{u,1}_\sigma)_H,u(\sigma)){|}d\sigma
%                       +\bigg{|}\int_{t_1}^{t_2}b(\sigma)X^{u,1}(\sigma-\tau)d\sigma\bigg{|}\\
%                       &&+ \int_{t_2}^{l}{|}F(\sigma,X^{u,1}(\sigma),
%                       (a,X^{u,1}_\sigma)_H,u(\sigma))-F(\sigma,X^{u,2}(\sigma),
%                       (a,X^{u,2}_\sigma)_H,u(\sigma)){|}d\sigma\\
%                &&+\sup_{s\in[t_2,l]}\bigg{|}\int_{t_2}^{s}b(\sigma)X^{u,1}(\sigma-\tau)-b(\sigma)X^{u,2}(\sigma-\tau)d\sigma\bigg{|}\\
                &\leq&
                |x_1(0)-x_2(0)|+L[1+(1+|a|_{H}\tau)\sup_{s\in[t_1,T]}|X^{1,u}(s)|+|a|_{H}|x_1|_H](t_2-t_1)\\
                &&+\sup_{s\in[0,T]}|b(s)|\bigg{|}\int^{t_2-t_1-\tau}_{-\tau}x_1(\theta)d\theta\bigg{|}
                        %+|b|_{W^{1,2}}\sup_{l\in[-\tau,0]}\bigg{|}\int^{l}_{-\tau}x_1(\sigma)d\sigma\bigg{|}(t_2-t_1)^{\frac{1}{2}}\\
%                &&
                +(1+\tau^{\frac{3}{2}}|a|_{W^{1,2}})\int^{l}_{t_2}\sup_{s\in[t_2,\sigma]}|X^{u,1}(s)-X^{u,2}(s)|d\sigma\\
                 && +L|a|_{W^{1,2}}\tau^{\frac{1}{2}}|T-t_2|\sup_{l\in[-\tau+t_2-t_1,0]}\bigg{|}\int^{0}_{l}x_{1}(\theta)-x_{2}(\theta+t_1-t_2)d\theta\bigg{|}\\
                &&+L|a|_{W^{1,2}}\tau^{\frac{1}{2}}|T-t_2|(t_2-t_1)^{\frac{1}{2}}[(t_2-t_1)^{\frac{1}{2}}\sup_{s\in[t_1,T]}|X^{u,1}(s)|+|x_2|]\\
                &&+\sup_{s\in[0,T]}|b(s)|\bigg{[}\bigg{|}\int^{0}_{t_2-t_1-\tau}x_1(\theta)-x_2(\theta+t_1-t_2)d\theta\bigg{|}
                  +(t_2-t_1)\sup_{s\in[t_1,T]}|X^{u,1}(s)|\\
                  &&~~+(t_2-t_1)^{\frac{1}{2}}|x_2|+\int^{l}_{t_2}\sup_{s\in[t_2,\sigma]}|X^{u,1}(s)-X^{u,2}(s)|d\sigma\bigg{]}.%+
%                16L^2[1+(2+\tau)\sup_{s\in[t_1,T]}|X^{1,u}(s)|^2+|x_1|^2_{H}](t_2-t_1)^2+16L^2(t_2-t_1)|x_1|_H^2\\
%                &&+16L^2(T-t_2)(T-t_2+2)\bigg{(}\int^{l}_{t_2}\sup_{s\in[t_1\vee
%                t_2,\sigma]}|X^{u,1}(s)-X^{u,2}(s)|^2d\sigma\\
%                &&+\int^{t_1-t_2}_{-\tau}|x_1(\sigma+t_2-t_1)-x_2(\sigma)|^2d\sigma+\int^{t_2}_{t_1}|X^{u,1}(\sigma)-x_2(\sigma-t_2)|^2d\sigma\bigg{)}.
\end{eqnarray*}
               Using the Gronwall-Bellman inequality,  we
                         obtain the following result, for a constant $C>0$,
\begin{eqnarray*}
                  && \sup_{s\in[t_1\vee t_2,T]}|X^{u,1}(s)-X^{u,2}(s)|\\
                  %&\leq&
%                  C(1+|x_1(0)|+|x_2(0)|+|x_1|_H+|x_2|_H){(}|x_1(0)-x_2(0)|+|t_2-t_1|+|t_2-t_1|^{\frac{1}{2}})
%                  \\
%                  &&+C\sup_{l\in[-\tau+t_2-t_1,0]}\bigg{|}\int^{0}_{l}x_{1}(\theta)-x_{2}(\theta+t_1-t_2)d\theta\bigg{|}\\
                  &\leq&
                  C(1+|x_1(0)|+|x_2(0)|+|x_1|_H+|x_2|_H){(}|x_1(0)-x_2(0)|+|t_2-t_1|+|t_2-t_1|^{\frac{1}{2}})
                  \\
                  &&+2C|x_1||t_2-t_1|^{\frac{1}{2}}+C\sup_{l\in[-\tau,0]}\bigg{|}\int^{0}_{l}x_{1}(\theta)-x_{2}(\theta)d\theta\bigg{|}.
\end{eqnarray*}
               Applying the supremum i.e.,  $\sup_{u\in {\cal{U}}[t_1\wedge
               t_2,T]}$, to both sides of the previous inequality, we obtain (2.6).
                    We can show that (2.7) holds using a similar (even simpler) procedure. \ \ $\Box$

%%%%%%%%%%%%%%%%%%%%%%%%%%%%%%%%%%%%%%%%%%%%%%%%%%%%%%%%%%%%%%%
\section{A DPP for optimal control problems}%Construction of random invariant manifold

\par
              In this section, we consider the controlled state
             equations:
$$
             X^u(s,t,x)=x(0)+\int_{t}^{s}F{(}\sigma,X^u(\sigma,t,x),(a,X^u_\sigma(t,x))_{H},u(\sigma){)}d\sigma
                    +\int_{t}^{s}b(\sigma)X^u(\sigma-\tau,t,x)d\sigma,\
             s\in [t,T], \eqno(3.1)
$$
       where $X^u_t=x\in {\cal{D}}$, and the cost function
$$
            ~J(t,x,u)=\int_{t}^{T}q(\sigma,X^u(\sigma,t,x),u(\sigma))d\sigma\\
                          +\phi(X^u(T,t,x)).\eqno(3.2)
$$
             Our purpose is to minimize the function $J$ over
             all controls  $u\in{\mathcal
                  {U}}[t,T]$.                           We define the function $V:[0,T]\times {\cal{D}}\rightarrow R$ by the following:
$$
                  V(t,x):=\inf_{u\in{\mathcal
                  {U}}[t,T]}J(t,x,u).\eqno(3.3)
$$
          The function $V$ is called the $value\ function$ of optimal
          control problem (3.1) and (3.2). The goal of this paper is
          to characterize this value function.
 %As in section 2, if $x\in{\cal{D}}$
%               and $x(0)=x^0$, we denote the cost function and the value function by $J(t,x,u)$ and $V(t,x)$, respectively.
 \par           We make the
              following assumptions.
\begin{hyp}
          \begin{description}
                  \item{(i)} The mappings $q:[0,T]\times R^d\times U\rightarrow R$ and
$\phi:R^d\rightarrow R$ are
                          measurable
                           and there exists a constant
                           $L>0$ , such that, for every $t\in [0,T], x\in R^d, u\in U$,
\begin{eqnarray*}
                          |q(t,x,u)|+|\phi(x)|\leq L(1+|x|).
\end{eqnarray*}
 \item{(ii)}            There exist a constant
                           $L>0$ and a local modulus of continuously $\rho$ such that, for every $t,s\in [0,T], x,y,\in R^d,   u\in U$,
\begin{eqnarray*}
                           |q(t,x,u)-q(s,y,u)|+|\phi(x)-\phi(y)|\leq L|x-y|+\rho(|s-t|,|x|\vee|y|).
\end{eqnarray*}
\end{description}
\end{hyp}
                Our first result is the local boundedness  and
                two kinds of
                continuities  of the value function.
\begin{theorem} Suppose that Hypothesis 2.1  and
                          Hypothesis 3.1 hold true. Then, there exists a constant $C_5>0$  such that,
                          for every $t,s\in [0,T]$, $x, y\in {\cal{D}}$,
   $$
                         |V(t,x)|\leq C_5(1+|x(0)|+|x|_H),\eqno(3.4)
$$
$$
                        |V(t,x)-V(t,y)|\leq
                        C_5{(}|x(0)-y(0)|+|x-y|_H{)},\eqno(3.5)
$$
    and
\begin{eqnarray*}
                        &&|V(t,x)-V(s,y)|\leq
                        C_5(1+|x(0)|+|y(0)|+|x|_H+|y|_H)\\
                        &&~~~~~~~~~~~~~~~~~~~~~~~~~~~~~~\times\bigg{(}|x(0)-y(0)|+|s-t|^{\frac{1}{2}}
                  +\sup_{l\in[-\tau,0]}\bigg{|}\int^{0}_{l}x(\theta)-y(\theta)d\theta\bigg{|}\bigg{)}.\
                  \ \ \ \ \ \ \ \
                  (3.6)
\end{eqnarray*}
\end{theorem}
 \par
{\bf  Proof}. \ \   %By the Hypothesis 3.1 (i), it is easy to obtain
%that (3.4) and (3.5) holds true.
                  We let $0\leq t\leq s\leq T$, $x,y\in {\cal{D}}$, by Hypothesis 3.1 (ii), (2.3) and (2.6), for any $u\in {\cal{U}}[t,T]$,
                  we have
\begin{eqnarray*}
                        &&|J(t,x,u)-J(s,y,u)|\\
                        %&\leq&|\phi(X^u(T,t,x))-\phi(X^u(T,s,y))|+\int_{t}^{s}q(\sigma,X^u(\sigma,t,x),u(\sigma))d\sigma\\
%                        &&+\int_{s}^{T}|
%                               q(\sigma,X^u(\sigma,t,x),u(\sigma)) -q(\sigma,X^u(\sigma,s,y),u(\sigma))|d\sigma\\
                        &\leq&(T+1)L\sup_{\sigma\in
                        [s,T]}[|X^u(\sigma,t,x)-X^u(\sigma,s,y)|+L\int_{t}^{s}(1+|X^u(\sigma,t,x)|)d\sigma\\
                        &\leq&
                        L(1+C_2(1+|x(0)|+|x|_H))(s-t)+(T+1)LC_3(1+|x(0)|+|y(0)|+|x|_H+|y|_H)\\
                        &&~~~\bigg{(}|x(0)-y(0)|+|s-t|^{\frac{1}{2}}
                  +\sup_{l\in[-\tau,0]}\bigg{|}\int^{0}_{l}x(\theta)-y(\theta)d\sigma\bigg{|}\bigg{)}.
\end{eqnarray*}
                  Thus, taking the infimum in $u\in {\cal{U}}[t,T]$, we
                  obtain (3.6). By the similar procedure, we can show
                    (3.4) and (3.5) hold true. The theorem is proved. \ \ $\Box$
\par
                    We note that $V(t,x)$ is not necessarily
                    Lipschitz continuous in $t$.
\par
                        Secondly, we present the following result, which is called  the dynamic programming
                        principle (DPP).
\begin{theorem}
                      Assume the Hypothesis 2.1   and Hypothesis
                      3.1  hold true. Then, for every $(t,x)\in [0,T)\times \cal{D}$ and $s\in [t,T]$, we have that
$$
             V(t,x)=\inf_{u\in {\cal{U}}[t,T]}\bigg{[}\int_{t}^{s}q(\sigma,X^u(\sigma,t,x),u(\sigma))d\sigma+V(s,X^u_s(t,x))\bigg{]}.  \eqno(3.7)
$$
\end{theorem}
 \par
{\bf  Proof}. \ \ First of all, for any  $u\in{\cal{U}}[s,T]$,
$s\in[t,T]$ and any  $u\in{\cal{U}}[t,s]$, by putting them
             concatenatively, we get $u\in{\cal{U}}[t,T]$.
                  Let us denote the right-hand side of (3.7) by
                  $\overline{V}(t,x)$. By (3.3), we have
\begin{eqnarray*}
           V(t,x)&\leq& J(t,x,{u})
              =\int_{t}^{s}q(\sigma,X^u(\sigma,t,x),u(\sigma))d\sigma+J(s,X^{{u}}_s(t,x),{{u}}),\
              u(\cdot)\in {\cal{U}}[t,T].
\end{eqnarray*}
             Thus, taking the infumum over $u(\cdot)\in
             {\cal{U}}[s,T]$, we obtain
\begin{eqnarray*}
           V(t,x)&\leq& \int_{t}^{s}q(\sigma,X^u(\sigma,t,x),u(\sigma))d\sigma+V(s,X^{{u}}_s(t,x)).
\end{eqnarray*}
             Consequently,
$$
                 V(t,x)\leq \overline{V}(t,x).
$$
            On the other hand, for any $\varepsilon>0$, there exists a $u^\varepsilon\in {\cal{U}}[t,T]$, such that
       \begin{eqnarray*}
           V(t,x)+\varepsilon&\geq& J(t,x,u^\varepsilon)\\
           &\geq&%\int_{t}^{s}q(\sigma,X^{u^\varepsilon}(\sigma,t,x),{u^\varepsilon}(\sigma))d\sigma
%           +\int_{s}^{T}q(\sigma,X^{u^\varepsilon}(\sigma,t,x),{u^\varepsilon}(\sigma))d\sigma
%                           +\phi(X^{u^\varepsilon}(T,t,x))\\
%              &=&
              \int_{t}^{s}q(\sigma,X^{u^\varepsilon}(\sigma,t,x),{u^\varepsilon}(\sigma))d\sigma
                           +J(s,X^{u^\varepsilon}_s(t,x),{u^\varepsilon})\\
              &\geq&\int_{t}^{s}q(\sigma,X^{u^\varepsilon}(\sigma,t,x),{u^\varepsilon}(\sigma))d\sigma
                           +V(s,X^{u^\varepsilon}_s(t,x))\geq\overline{V}(t,x).
\end{eqnarray*}
             Hence, (3.7) follows.\ \ $\Box$
\par
                  Our next goal is to derive the so-called
                  $Hamilton-Jacobi-Bellman\
                  equation$ for the value function $V$. To begin
                  with, let us introduce the operator ${\cal{S}}$.
                     For a  Borel
                      measurable function $f:{\cal{D}}\rightarrow
                      R$, we define
$$
                       {\mathcal
                       {S}}(f)(x)=\lim_{h\rightarrow0^+}\frac{1}{h}[f(\widehat{{x}}_h)-f(x)],\
                       \ \ x\in {\cal{D}},
$$
                     where $\widehat{{x}}:[-\tau,T]\rightarrow R^d$
                     is an extension of $x$ defined by
\begin{equation*}
 \widehat{{x}}(s)=
\begin{cases}
x(s),  \ \ \ \ s\in
                                                         [-\tau,0),\\
x(0), \ \  \ \  s\geq 0,
\end{cases}
\end{equation*}
%\begin{eqnarray*}
%                          \widehat{{x}}(s)=\cases{ x(s),  \ \ \ \ s\in
%                                                         [-\tau,0),\cr
%                                            x(0), \ \  \ \  s\geq 0,\cr}
%\end{eqnarray*}
                             and $ \widehat{{x}}_s$ is defined by
$$
                             \widehat{{x}}_s(\theta)=\widehat{{x}}(s+\theta),\
                              \ \ \theta\in[-\tau,0].
$$
                        We denote by $\widehat{D}(\mathcal
                      {S})$ the domain of the operator ${\mathcal
                      {S}}$, be the set of $f:{\cal{D}}\rightarrow
                      R$ such that the above limit exists for all $x \in
                      {\cal{D}}$. Define $D({\mathcal {S}})$ as
                      the set of all functions $g: [0,T]\times {\cal{D}}\rightarrow
                      R$ such that $g(t,\cdot)\in \widehat{D}(\mathcal
                      {S})$ for all $t\in [0,T]$.  For simplicity, we define
\begin{eqnarray*}
\Phi&=&\{\varphi\in C^{1} ([0,T]\times  {\cal{D}})\cap
             D({\cal{S}})| \ \exists\ \varphi_0\in C^{1}([0,T]\times R^d\times H),\\
                             &&~ \ \mbox{such that }  \varphi(t,x)=\varphi_0(t,x(0),x),  \ \ \ \forall (t,x)\in[0,T]\times {\cal{D}}
\}.
\end{eqnarray*}
\begin{theorem}
                      Let $V$ denote the value function defined by (3.3), if the function $V(t,x)\in \Phi$. Then,
                 $V(t,x)$ satisfies the following HJB equation:
$$
\begin{cases}
\frac{\partial}{\partial t} V(t,x)+{\mathcal
                        {S}}(V)(t,x)+ H(t,x,\nabla_x V(t,x))= 0, \ \  t\in
                               [0,T],\ \ x\in {\mathcal{D}}, \\
 V(T,x)=\phi(x(0)),
\end{cases}\eqno(3.8)
$$
%$$
%          \cases{\frac{\partial}{\partial t} V(t,x)+{\mathcal
%                        {S}}(V)(t,x)+ H(t,x,\nabla_x V(t,x))= 0, \ \  t\in
%                               [0,T],\ \ x\in {\mathcal{D}}, \cr
%          V(T,x)=\phi(x(0)),\cr} \eqno(3.8)
%$$
                where
\begin{eqnarray*}
                                H(t,x,p)=\inf_{u\in{
                                         {U}}}[\langle
                        [F(t,x(0),(a,x)_{H},u)+b(t)x(-\tau)]1_0(t),p\rangle+q(t,x(0),u)],\  (t,x,p)\in [0,T]\times {\cal{D}}\times
                        {\cal{D}}^*.
\end{eqnarray*}
                 Here the function $1_{0}:[-\tau,0]\rightarrow R$
                       is the character function of $\{0\}$.
\end{theorem}
%
%                      We define ${\mathcal {L}}_t$ by
%$$
%                        {\mathcal {L}}_t[\phi](x)={\mathcal
%                        {S}}(\phi)(x)+\langle
%                        F(t,x(0),x(-\tau),x,u),{\nabla_{0}\phi(x)}\rangle,
%$$
%                        where  $\nabla_{x_0}\phi$, is
%                        first derivatives
%                         of $\phi$ with respect to $x_0$.
%\par\vskip2mm
                In order to prove this theorem we need the following
                lemma.
\begin{lemma}
                          Suppose that Hypothesis 2.1 holds.  If $g\in \Phi$, then, for each
                       $(t,x)\in [0,T)\times{\mathcal{D}}$, the following convergence holds uniformly in $u(\cdot)\in {\cal{U}}[t,T]$:
$$
                       \lim_{\epsilon\rightarrow0^+}
                       \bigg{[}\frac{g(t+\epsilon,X^u_{t+\epsilon})-g(t,x)}{\epsilon}-g_t(t,x)-{\mathcal
                        {S}}(g)(t,x)-\langle{\nabla_{x}g(t,x)},\overline{W}(t,x,\epsilon)1_0(t)\rangle\bigg{]}=0,
                     \eqno(3.9)
$$
where we let  $\overline{W}(t,x,\epsilon)$ denote
$\frac{1}{\epsilon}\int^{t+\epsilon}_{t}
                               F{(}t,x(0),(a,x)_H,u(\sigma){)}d\sigma
                                 +b(t)x(-\tau)$.
\end{lemma}
                     %Here $g_0\in C^{1}([0,T]\times R^d\times H)$ \ such that \  $g(t,x)=g_0(t,x(0),x)$,  for all $(t,x)\in[0,T]\times
%                     {\cal{D}}$, where  $\nabla_{x^0}g_0$ denotes first derivative with respect to the second variable of $g_0$.
\par
   {\bf  Proof}. \ \
                      Since $g\in C^{1}([0,T]\times
                 {\mathcal{D}})$,  by Taylor's theorem we get
                 that
\begin{eqnarray*}
                    g(t+\epsilon,X^u_{t+\epsilon})-g(t,x)&=&g(t+\epsilon,X^u_{t+\epsilon})-g(t,X^u_{t+\epsilon})
                           +g(t,\widehat{x}_{t+\epsilon})-g(t,x)\\
                           &&+\langle\nabla_x g(t,\widehat{x}_{t+\epsilon}),
                            X^u_{t+\epsilon}-\widehat{x}_{t+\epsilon})\rangle+o(|X^u_{t+\epsilon}-\widehat{x}_{t+\epsilon}|), \
                            t\in [0,T-\epsilon]. \ (3.10)
\end{eqnarray*}
 Again by $g\in C^{1}([0,T]\times
                      {\mathcal{D}})$, we have that
$$
                        \lim_{\epsilon\rightarrow0^+}\sup_{u(\cdot)\in {\cal{U}}[t,T]}\bigg{|}\frac{1}{\epsilon}[g(t+\epsilon,X^u_{t+\epsilon})
                                    -g(t,X^u_{t+\epsilon})] -\frac{\partial}{\partial t}g(t,x)\bigg{|}=0.\eqno(3.11)
$$
             From $g\in D({\cal{S}})$, it follows that
$$
                         \lim_{\epsilon\rightarrow0^+}\frac{1}{\epsilon}(g(t,\widehat{x}_{t+\epsilon})-g(t,x))={\cal{S}}(g)(t,x).\eqno(3.12)
$$
    By the definitions of $X^u_s$ and $\widehat{x}_s$, we have
                        that, for every $\epsilon\in[0,T-t]$,
\begin{equation*}
 X^u_{t+\epsilon}(\theta)-\widehat{x}_{t+\epsilon}(\theta)=
\begin{cases}
\int^{{t+\epsilon}+\theta}_{t}F(\sigma,X^u(\sigma),(a,X^u_\sigma)_H,u(\sigma))+b(\sigma)X^u(\sigma-\tau)d\sigma,
                        \ \
                                      {\epsilon}+\theta\geq 0, \\
0,         \ \ \ \ \ \ \ \ \ \ \ \ \ \ \ \ \ \ \ \ \ \ \ \ \ \ \ \ \ \ \ \  \ \ \ \ \ \ \ \  \ \ \ \ \ \ \ \
                   \ \ \ \ \ \ \ \ \ \ \ \ \ \ \ \
                   \ \ \  \ \ \ \ \
                     {\epsilon}+\theta< 0.
\end{cases}
\end{equation*}
%$$
%                        X^u_{t+\epsilon}(\theta)-\widehat{x}_{t+\epsilon}(\theta)=
%                        \cases{\int^{{t+\epsilon}+\theta}_{t}F(\sigma,X^u(\sigma),(a,X^u_\sigma)_H,u(\sigma))+b(\sigma)X^u(\sigma-\tau)d\sigma,
%                        \ \
%                                      {\epsilon}+\theta\geq 0, \cr
%                   0,         \ \ \ \ \ \ \ \ \ \ \ \ \ \ \ \ \ \ \ \ \ \ \ \ \ \ \ \ \ \ \ \  \ \ \ \ \ \ \ \  \ \ \ \ \ \ \ \
%                   \ \ \ \ \ \ \ \ \ \ \ \ \ \ \ \
%                   \ \ \  \ \ \ \ \
%                     {\epsilon}+\theta< 0.\cr}
%$$
                      Thus, the following convergence holds uniformly in $u(\cdot)\in {\cal{U}}[t,T]$:
\begin{equation*}
\lim_{\epsilon\rightarrow 0^+}[\frac{1}{\epsilon}(X^u_{t+\epsilon}(t,x)-\widehat{x}_{t+\epsilon})](\theta)=
\begin{cases}
\lim_{\epsilon\rightarrow 0^+}\frac{1}{\epsilon}\int^{{t+\epsilon}}_{t}
                                       F(\sigma,X^u(\sigma),(a,X^u_\sigma)_H,u(\sigma))\\
 ~~~~~~~~~~~~~~~~~~~+b(\sigma)X^u(\sigma-\tau)d\sigma,\
                                \ \ \ \ \ \ \ \ \ \ \theta=0, \\
         0,         \ \ \ \ \ \ \ \ \ \ \ \ \ \ \  \ \ \ \ \ \ \ \ \ \
                   \ \ \ \ \ \ \ \ \ \ \  \ \ \ \ \ \ \ \ \ \ \
                     -\tau\leq\theta< 0.                        
\end{cases}
\end{equation*}
%$$
%                        \lim_{\epsilon\rightarrow 0^+}[\frac{1}{\epsilon}(X^u_{t+\epsilon}(t,x)-\widehat{x}_{t+\epsilon})](\theta)=
%                        \cases{\lim_{\epsilon\rightarrow 0^+}\frac{1}{\epsilon}\int^{{t+\epsilon}}_{t}
%                                       F(\sigma,X^u(\sigma),(a,X^u_\sigma)_H,u(\sigma))\cr
%                                       ~~~~~~~~~~~~~~~~~~~+b(\sigma)X^u(\sigma-\tau)d\sigma,\
%                                \ \ \ \ \ \ \ \ \ \ \theta=0, \cr
%                   0,         \ \ \ \ \ \ \ \ \ \ \ \ \ \ \  \ \ \ \ \ \ \ \ \ \
%                   \ \ \ \ \ \ \ \ \ \ \  \ \ \ \ \ \ \ \ \ \ \
%                     -\tau\leq\theta< 0.\cr}
%$$
                   By (2.3), we
                     get that, there exists a constant $C>0$ independent of $u(\cdot)$ such that
\begin{eqnarray*}
                     &&\sup_{\theta\in [-\tau,0]}|[\frac{1}{\epsilon}(X^u_{t+\epsilon}-\widehat{x}_{t+\epsilon})](\theta)|
                     \leq C.
\end{eqnarray*}%
%                           For every $t\in[0,T]$, let
%                           $\nabla_{(x^0,x)}g_0(t,\cdot)$ denote the
%                           first Fr\'{e}chet derivative of $g_0$
%                           with respect to $(x^0,x)\in R^d\times H$.
%                           It is clearly that
%                           $\nabla_{(x^0,x)}g_0(t,\cdot)$ is
%                           continuous.
                           Thus, by the continuity of
                           $\nabla_{x^0}g_0(t,\cdot,x)$ and
                           $\nabla_{x}g_0(t,x^0,\cdot)$,
                            we obtain
\begin{eqnarray*}
                             &&\lim_{\epsilon\rightarrow0^+}\sup_{u(\cdot)\in {\cal{U}}[t,T]}\bigg{|}\frac{1}{\epsilon}\langle\nabla_{x} g(t,\widehat{x}_{t+\epsilon}),
                                         X^u_{t+\epsilon}-\widehat{x}_{t+\epsilon}\rangle-\frac{1}{\epsilon}\langle\nabla_{x}g(t,{x}),
                              X^u_{t+\epsilon}-\widehat{x}_{t+\epsilon}\rangle\bigg{|}\\
                             %&\leq&C\lim_{\epsilon\rightarrow0^+}|\nabla_{x}g(t,\widehat{x}_{t+\epsilon})-\nabla_{x}g(t,{x})|\\
                             &\leq&C\lim_{\epsilon\rightarrow0^+}|\nabla_{x^0}g_0(t,x(0),\widehat{x}_{t+\epsilon})
                                               +\nabla_{x}g_0(t,x(0),\widehat{x}_{t+\epsilon})-\nabla_{x^0}g_0(t,x(0),{x})-\nabla_{x}g_0(t,x(0),{x})|\\
                             &=&0.%(\mbox{there is a bug}, \widehat{x}_{t+\epsilon}\nrightarrow x \ \mbox{ as}\  \varepsilon\rightarrow0)
\end{eqnarray*}
                          Therefore, by the above inequality, % by Hypothesis 2.1 and Theorem 2.3,
                          we obtain that
\begin{eqnarray*}
                       ~~~~&&\lim_{\epsilon\rightarrow0^+}\sup_{u(\cdot)\in
                       {\cal{U}}[t,T]}\bigg{|}\frac{1}{\epsilon}\langle\nabla_{x}g(t,\widehat{x}_{t+\epsilon}),
                               X^u_{t+\epsilon}-\widehat{x}_{t+\epsilon}\rangle-\langle\nabla_{x}g(t,{x}),\overline{W}(t,x,\epsilon)1_0(t)\rangle \bigg{|}\\
                       &=&\lim_{\epsilon\rightarrow0^+}\sup_{u(\cdot)\in
                       {\cal{U}}[t,T]}\bigg{|}\frac{1}{\epsilon}\langle\nabla_{x}g(t,x),
                               X^u_{t+\epsilon}-\widehat{x}_{t+\epsilon}\rangle-\langle\nabla_{x}g(t,{x}),\overline{W}(t,x,\epsilon)1_0(t)\rangle \bigg{|}\\
                        &\leq&\lim_{\epsilon\rightarrow0^+}\sup_{u(\cdot)\in
                       {\cal{U}}[t,T]}\bigg{|}\langle\nabla_{x}g_0(t,x(0),x),
                               \frac{1}{\epsilon}(X^u_{t+\epsilon}-\widehat{x}_{t+\epsilon})
                               \rangle \bigg{|}+\lim_{\epsilon\rightarrow0^+}\sup_{u(\cdot)\in
                       {\cal{U}}[t,T]}\bigg{|}\langle\nabla_{x^0}g_0(t,x(0),x),\\
                       &&~~~~~~~~~~~~
                               \frac{1}{\epsilon}\int^{t+\epsilon}_{t}[F{(}\sigma,X^u(\sigma),(a,X^u_\sigma)_H,u(\sigma){)}
                               + b(\sigma)X^u(\sigma-\tau)]d\sigma-
                                                   \overline{W}(t,x,\epsilon)
                               \rangle \bigg{|}\\
                       %&=&\lim_{\epsilon\rightarrow0^+}\sup_{u(\cdot)\in
%                       {\cal{U}}}\bigg{|}\langle\nabla_{x}g(t,{x}),
%                               1_0\frac{1}{\epsilon}\int^{{t+\varepsilon}}_{t}
%                                       F(\sigma,X^u(\sigma),X^u(\sigma-\tau),X^u_\sigma,u(\sigma))d\sigma
%                                      \\
%                       &&~~~~~~~~~~~~~~~~~~~~~~~~-1_0\frac{1}{\epsilon}\int^{t+\epsilon}_{t}F(t,x(0),x(-\tau),x,u(\sigma))d\sigma\rangle\bigg{|}\\
%                        &\leq&L\lim_{\epsilon\rightarrow0^+}\sup_{u(\cdot)\in {\cal{U}}}{|}\nabla_{x}g(t,{x})|
%                               \int^{{t+\varepsilon}}_{t}\frac{|X^u(\sigma)-x(0)|+|X^u(\sigma-\tau)-x(-\tau)|+|X^u_\sigma-x|}{\epsilon}d\sigma\\
%                               &&+\lim_{\epsilon\rightarrow0^+}\sup_{u(\cdot)\in
%                               {\cal{U}}}{|}\nabla_{x}g(t,{x})|
%                               \int^{{t+\varepsilon}}_{t}\frac{\gamma(|\sigma-t|,|x(0)|+|x(-\tau)|+|x|)}{\epsilon}d\sigma
                       & =&0.\ \  \ \ \ \ \ \ \ \ \ \ \ \ \ \ \ \  \ \ \ \ \ \  \ \ \ \ \ \ \ \ \ \ \ \ \ \ \ \  \ \ \ \ \
                       \  \ \ \ \ \ \ \ \ \ \ \ \ \ \ \ \  \ \ \ \ \
                       \  \ \ \ \ \ \ \ \ \ \ \ \ \ \ \  \ \ \ \ \ \  \ \ \  \ \ \ \ \ \ \ \ \  \ \ \ \ \    (3.13)
\end{eqnarray*}
                Hence, dividing by $\varepsilon$ in (3.10) and
                      sending  $\epsilon\rightarrow0^+$, and putting together the results
                          of (3.11), (3.12) and (3.13), we
                          finally obtain (3.9).\ \ $\Box$
\par
              From the above lemma,  the following two lemmas hold
              true, which will be used in the proof of uniqueness
              result for viscosity solution.
\begin{lemma}
                          Suppose that Hypothesis 2.1 holds. If $g(t,x)=g_0(t,|x|^2_H),\  (t,x)\in [0,T)\times  {\cal{D}}$,
                           where $g_0\in C^{1}([0,T]\times R)$.
                           Then the following  holds:
$$
                       \frac{g(t+\epsilon,X^u_{t+\epsilon})
                                      -g(t,x)}{\epsilon}=
                                      g_t(t,x)+g'_0(t,|x|_H^2)(x^2(0)-x^2(-\tau))+o(1),\eqno(3.14)
$$
                   where $o(1)$ is uniformly in $u(\cdot)\in {\cal{U}}[t,T]$.
\end{lemma}
\par
   {\bf  Proof}. \ \                  By Lemma 3.5, we only need to show that %$g\in \Phi$ and   By Taylor's Theorem we get
%                 that
%\begin{eqnarray*}
%                    g(t+\epsilon,X^u_{t+\epsilon})-g(t,x)
%                    &=&g(t+\epsilon,X^u_{t+\epsilon})-g(t,X^u_{t+\epsilon})
%                           +g(t,\widehat{x}_{t+\varepsilon})-g(t,x)\\
%                    &&+\nabla_{x} g(t,\widehat{x}_{t+\varepsilon})(X^u_{t+\epsilon}-\widehat{x}_{t+\varepsilon})
%                     +o(|X^u_{t+\epsilon}-\widehat{x}_{t+\varepsilon}|), \
%                            t\in[0,T-\epsilon].
%\end{eqnarray*}
%                        By the similar   procedure of Lemma 3.4, we get that
$$
                        \lim_{\epsilon\rightarrow0^+}\sup_{u(\cdot)\in {\cal{U}}[t,T]}\bigg{|}\frac{1}{\epsilon}[g(t+\epsilon,X^u_{t+\epsilon})
                                         -g(t,X^u_{t+\epsilon})]
                                         -g_t(t,x)\bigg{|}=0,\eqno(3.15)
$$
$$
   \lim_{\epsilon\rightarrow0^+}\frac{1}{\epsilon}\sup_{u(\cdot)\in {\cal{U}}[t,T]}|\nabla_xg(t,\widehat{x}_{t+\epsilon})
                                   (X^u_{t+\epsilon}-\widehat{x}_{t+\epsilon})|=0,\eqno(3.16)
$$
and
$$
             \lim_{\epsilon\rightarrow 0}\frac{1}{\epsilon}(g(t,\widehat{x}_{t+\epsilon})-g(t,x))=
             g'_0(t,|x|_H^2)(x^2(0)-x^2(-\tau)).\eqno(3.17)
$$
 By the similar procedure of Lemma 3.5, we get (3.15) and (3.16)
 hold true.
                      Now, let us prove (3.17).
             From the definition of $g$, it follows that
 \begin{eqnarray*}
                        && \frac{1}{\epsilon}(g(t,\widehat{x}_{t+\epsilon})-g(t,x))
                               = \frac{1}{\epsilon}(g_0(t,|\widehat{x}_{t+\epsilon}|_H^2)-g_0(t,|x|_H^2))\\
                        &=&\frac{1}{\epsilon}\int^{1}_{0}g'_0(t,|x|_H^2+s(|\widehat{x}_{t+\epsilon}|_H^2-|x|_H^2))
                             (|\widehat{x}_{t+\epsilon}|_H^2-|x|_H^2)ds\\
                        &=& \frac{1}{\epsilon}\int_{-\tau}^{0}[\widehat{x}^2_{t+\epsilon}(\theta)-x^2(\theta)]d\theta
                        \int^{1}_{0}g'_0(t,|x|_H^2+s(|\widehat{x}_{t+\epsilon}|_H^2-|x|_H^2))ds\\
                        &=&
                        \bigg{(}{x}^2(0)-\frac{1}{\epsilon}\int^{-\tau+\epsilon}_{-\tau}x^2(\theta)d\theta\bigg{)}
                        \int^{1}_{0}g'_0(t,|x|_H^2+s(|\widehat{x}_{t+\epsilon}|_H^2-|x|_H^2))ds.
\end{eqnarray*}
             Letting $\epsilon\rightarrow 0$, we obtain (3.17).\ \ $\Box$
%\par
%We define the $|\cdot|_{B}$-norm on ${\cal{D}}$ by
%$$
%                              |x|^2_{B}:=\int^{0}_{-\tau}(Bx)^2(s)ds,
%$$
%              where
% $$
%            ( Bx)(s)=\int^{0}_{s}x(\theta)d\theta,\ \ s\in [-\tau,0].
% $$
\begin{lemma}
                          Suppose that Hypothesis 2.1 holds. If $\psi(x)=\psi_0(|x-\hat{a}|_B^2), \hat{a},x\in {\cal{D}}$, where $\psi_0\in C^1(R)$. Then
                             the following  holds:
\begin{eqnarray*}
                         \frac{1}{\epsilon}\psi(X^u_{t+\epsilon})-\psi(x)
                         %&=&2\psi'_0(|(x-a)|_B^2)\bigg{(}-\int^{0}_{-\tau}x(s)\int^{0}_{s}(x(\theta)-a(\theta))d\theta
%                                 ds+x_0\int^{0}_{-\tau}\int^{0}_{s}(x(\theta)-a(\theta))d\theta ds\bigg{)}+o(1)\\
                         =2\psi'_0(|(x-\hat{a})|_B^2)( B(x-\hat{a}),x(0)1_{[-\tau,0]}-x)_H+o(1),
\end{eqnarray*}
                   where $o(1)$ is uniformly in $u(\cdot)\in
                   {\cal{U}}[t,T]$.  Here the function $1_{[-\tau,0]}$
                       is the character function of $[-\tau,0]$.
\end{lemma}
\par
   {\bf  Proof}. \ \  By Lemma 3.5, we only need to show that
%$$
%                    \frac{1}{\epsilon}[\psi(X^u_{t+\epsilon})-\psi(x)]
%                    = \frac{1}{\epsilon}(\psi(\widehat{x}_{t+\varepsilon})-\psi(x))+\frac{1}{\epsilon}\nabla_x \psi(\widehat{x}_{t+\varepsilon})(X^u_{t+\epsilon}
%                                -\widehat{x}_{t+\varepsilon})
%                       +o(|X^u_{t+\epsilon}-\widehat{x}_{t+\varepsilon}|),
%                            t\in[0,T-\epsilon],
%$$
%              and
$$
          \lim_{\epsilon\rightarrow0^+}\frac{1}{\epsilon}\sup_{u(\cdot)\in {\cal{U}}[t,T]}|\nabla_x\psi(\widehat{x}_{t+\epsilon})
                                   (X^u_{t+\epsilon}-\widehat{x}_{t+\epsilon})|=0,\eqno(3.18)
$$
and
$$
                        \lim_{\epsilon\rightarrow0} \frac{1}{\epsilon}(\psi(\widehat{x}_{t+\epsilon})-\psi(x))
                                =2\psi'_0(|(x-\hat{a})|_B^2)(
                                B(x-\hat{a}),x(0)1_{[-\tau,0]}-x)_H.\eqno(3.19)
$$
    By the similar procedure of Lemma 3.5, we can obtain (3.18).
                     Now let us show (3.19) hold true.  By the definition of $\psi$, we have that, for some
                         $s\in (0,1)$,
 \begin{eqnarray*}
                       ~~~~~~~~~~~~ && \frac{1}{\epsilon}(\psi(\widehat{x}_{t+\epsilon})-\psi(x))
                                =\frac{1}{\epsilon}[\psi'(x+s(\widehat{x}_{t+\epsilon}-x))(\widehat{x}_{t+\epsilon}-x)]\\
                        &=& \frac{2}{\epsilon}\psi'_0(|(x+s(\widehat{x}_{t+\epsilon}-x)-\hat{a})|_B^2)
                             \langle
                             B(x+s(\widehat{x}_{t+\epsilon}-x)-\hat{a}),B(\widehat{x}_{t+\epsilon}-x)\rangle.\
                             \ \ \ \ \ \ \ \ \  \ (3.20)
\end{eqnarray*}
            On the other hand, we have that
\begin{eqnarray*}
                       &&\frac{1}{\epsilon}\langle B(x-\hat{a}),B(\widehat{x}_{t+\epsilon}-x)\rangle
                           =\frac{1}{\epsilon}\int^{0}_{-\tau}\int^{0}_{s}(x(\theta)-\hat{a}(\theta))d\theta
                                    \int^{0}_{s}(\widehat{x}_{t+\epsilon}(\theta)-x(\theta))d\theta ds\\
                      % &=&\frac{1}{\epsilon}\int^{0}_{-\tau}\int^{0}_{s}(x(\theta)-\hat{a}(\theta))d\theta
%                                  \bigg{(}  \int^{s\vee(-\epsilon)}_{s}x(\theta+\epsilon)d\theta+\int_{s\vee(-\epsilon)}^{0}x(0)d\theta
%                                  -\int^{0}_{s}x(\theta)d\theta \bigg{)}ds\\
                       &=&\frac{1}{\epsilon}\int^{0}_{-\tau}\int^{0}_{s}(x(\theta)-\hat{a}(\theta))d\theta
                                  \bigg{(}  -\int^{(s+\epsilon)\wedge0}_{s}x(\theta)d\theta+((-s)\wedge\epsilon)
                                  x(0)
                                  \bigg{)}ds\\
                        &\rightarrow&\int^{0}_{-\tau}(x(0)-x(s))\int^{0}_{s}(x(\theta)-\hat{a}(\theta))d\theta
                                 ds =( B(x-\hat{a}),x(0)1_{[-\tau,0]}-x)_H\ \ \ \ \mbox{as}\ \
                                 \epsilon\rightarrow0,
\end{eqnarray*}
                    and
\begin{eqnarray*}
                            &&\frac{1}{\epsilon}\langle B((\widehat{x}_{t+\epsilon}-x)),B(\widehat{x}_{t+\epsilon}-x)\rangle
                            =\frac{1}{\epsilon}\int^{0}_{-\tau}\bigg{(}\int^{0}_{s}\widehat{x}_{t+\epsilon}(\theta)-x(\theta)d\theta\bigg{)}^2ds\\
                            &=&\frac{1}{\epsilon}\int^{0}_{-\tau}\bigg{(}-\int^{(s+\epsilon)\wedge0}_{s}x(\theta)d\theta+((-s)\wedge\epsilon) x(0)\bigg{)}^2ds\\
                            &\leq&\frac{2}{\epsilon}\int^{0}_{-\tau}\bigg{(}\epsilon^2x^2(0)+\epsilon\int^{(s+\epsilon)\wedge0}_{s}x^2(\theta)d\theta \bigg{)}ds
                            \rightarrow0\ \ \ \ \mbox{as}\ \ \epsilon\rightarrow0.
 \end{eqnarray*}
  Letting $\epsilon\rightarrow0$ in (3.20), we get (3.19).\ \  $\Box$
  \begin{remark}
                \rm{ We note that $\widehat{D}(\mathcal
                      {S})$ and $D({\mathcal {S}})$ are not empty.
                      In fact, by (3.17) and (3.19), we have $g_0(|x|_H^2),\  g_0(|x-a|_B^2)\in\widehat{D}(\mathcal
                      {S})$, if $g_0\in C^1(R)$ and $ a\in {\cal{D}}$.
                      Moreover, $g_0(|x|_H^2)+l(t)$ and
                      $g_0(|x-a|_B^2)+l(t)$ belong to $D({\mathcal
                      {S}})$, if $l\in C^1(R)$.}
  \end{remark}
 \par
{\bf  Proof of Theorem 3.4}. \ \
                  First of all, by the  definition of $V$, we have that $V(T,x)=\phi(x(0))$. Next, we fix a $u\in U$ and $x\in {\cal{D}}$,
                  from (3.7),  it follows that
\begin{eqnarray*}
                 0\leq \int_{t}^{s}q(\sigma,X^u(\sigma,t,x),u)d\sigma
                                +V(s,X^u_s(t,x))
                 -V(t,x).
\end{eqnarray*}
              By Lemma 3.5, the above inequality implies that
\begin{eqnarray*}
    0&\leq&\lim_{s\rightarrow t^+}\frac{1}{s-t}\bigg{[}\int_{t}^{s}q(\sigma,X^u(\sigma,t,x),u)d\sigma+V(s,X^u_s(t,x))-V(t,x)\bigg{]}\\
                       &=&\frac{\partial}{\partial t} V(t,x)+{\mathcal
                        {S}}(V)(t,x)+ \langle
                        [F(t,x(0),(a,x)_{H},u)+b(t)x(-\tau)]1_0(t),\nabla_xV(t,x)\rangle+q(t,x(0),u).
\end{eqnarray*}
      Thus, we have that
$$
    0\leq\frac{\partial}{\partial t} V(t,x)+{\mathcal
                        {S}}(V)(t,x)+ H(t,x,\nabla_x V(t,x)).\eqno(3.21)
$$
          On the other hand,  let $x\in {\cal{D}}$ be fixed. For any $\varepsilon>0$ and $s>t$, by (3.7), there exists a
                         $\tilde{u}\equiv u^{\varepsilon,s}\in {\cal{U}}[t,T]$ such that
\begin{eqnarray*}
    \varepsilon(s-t)&\geq&\int_{t}^{s}q(\sigma,X^{\tilde{u}}(\sigma,t,x),\tilde{u}(\sigma))d\sigma+V(s,X^{\tilde{u}}_s(t,x))-V(t,x)\\
                     &=& \frac{\partial }{\partial t}V(t,x)(s-t)+{\mathcal
                        {S}}(V)(t,x)(s-t)+\int_{t}^{s}q(t,x(0),\tilde{u}(\sigma))d\sigma\\
                        &&+\langle{\nabla_{x}V(t,x)},\int^{s}_{t}
                        [F(t,x(0),(a,x)_{H},{\tilde{u}}(\sigma))+b(t)x(-\tau)]d\sigma1_0(t)\rangle+o(|s-t|)\\
                     &\geq& \frac{\partial }{\partial t}V(t,x)(s-t)+{\mathcal
                        {S}}(V)(t,x)(s-t)+H(t,x,\nabla_x V(t,x))(s-t)+o(|s-t|).
\end{eqnarray*}
Then, dividing through by $s-t$ and letting $s-t\rightarrow0$, we
have that
\begin{eqnarray*}
    \varepsilon&\geq&\frac{\partial}{\partial t} V(t,x)+{\mathcal
                        {S}}(V)(t,x)+ H(t,x,\nabla_x V(t,x)).
\end{eqnarray*}
               Combining with (3.21), we get the desired result.\ \ $\Box$

%%%%%%%%%%%%%%%%%%%%%%%%%%%%%%%%%%%%%%%%%%%%%%%%%%%%%%%%%
\section{Viscosity solution of HJB equations: Existence theorem.}
\par
                                   In this section, we are going to introduce the notion of viscosity
                                   solution.   M. G. Crandall and P. L. Lions \cite{cra3}, \cite{cra4}, \cite{cra5}, \cite{cra6}, \cite{cra7}
                                   systematically introduced  the basic theories of viscosity
              solutions for HJB equations in infinite dimensions. The proof of the uniqueness is mainly based on the weak compactness of separable
              Hilbert space (see also \cite{li}). We note that the  HJB
              equation (3.8) is defined in ${\cal{D}}$, which doesn't have the weak
              compactness. Thus, we need to give a new notion of
              viscosity solution of (3.8). To begin with, let us introduce  the following key lemma that will be used
              in the proof of the uniqueness of viscosity solutions.
\begin{lemma} {\bf(Left maximization principle)} Let $Q$ be a bounded open subset of $R^d$
and
   let $v: [0,T]\times {\cal{D}}\times[0,T]\times {\cal{D}}\rightarrow R$ be continuous, and
   there exists an integer $k>0$ such that, for every
   $t,s\in[0,T],x,x_1,x_2,y,y_1,y_2\in{\cal{D}}$,
$$
             |v(t,x,s,y)|\leq L(1+|x|_{H}+|x(0)|+|y|_{H}+|y(0)|)^k,
$$
$$
             |v(t,x_1,s,y_1)-|v(t,x_2,s,y_2)|\leq
             L(|x_1-x_2|_{H}+|x_1(0)-x_2(0)|+|y_1-y_2|_{H}+|y_1(0)-y_2(0)|)^k.\eqno(4.1)
$$ Then,
   for each $(t_0,x_0,s_0,y_0)\in [0,T]\times {\cal{D}}_{\bar{Q}}\times [0,T]\times {\cal{D}}_{\bar{Q}}$,
   there exists $(\bar{t},\bar{x},\bar{s},\bar{y})\in [t_0,T]\times
   {\cal{D}}_{\bar{Q}}\times [s_0,T]\times
   {\cal{D}}_{\bar{Q}}$, such that
   $(\bar{t},\bar{x},\bar{s},\bar{y})=(t_0,x_0,s_0,y_0)\otimes
   (\bar{t},\bar{x},\bar{s},\bar{y})$,  %$(\bar{s},\bar{y})=(s_0,y_0)\otimes
%   (\bar{s},\bar{y})$,
  $u(\bar{t},\bar{x},\bar{s},\bar{y})\geq u(t_0,x_0,s_0,y_0)$, and
   $$
                         v(\bar{t},\bar{x},\bar{s},\bar{y})=\sup_{(t,x,s,y)\in [\bar{t},T]\times
                         {\cal{D}}_{\bar{Q}}\times [\bar{s},T]\times
                         {\cal{D}}_{\bar{Q}}}
                         v((\bar{t},\bar{x},\bar{s},\bar{y})\otimes
                         ({t},{x},s,y)).\eqno(4.2)
   $$
\end{lemma}
 \par
{\bf  Proof}. \ \  Without loss of generality, we can assume that
                   $v(t_0, x_0,s_0,y_0)\geq v(t_0, x_0+e1_0(\cdot),s_0, y_0+l1_0(\cdot))$ for all $e,l\in R^d$ such
                   that $e+x_0(0),l+y_0(0)\in \bar{Q}$. We set $m_0=v(t_0,x_0,s_0,y_0)$ and
$$
                               \bar{m}_0:=\sup_{(t,x,s,y)\in [t_0,T]\times
                         {\cal{D}}_{\bar{Q}}\times [s_0,T]\times
                         {\cal{D}}_{\bar{Q}}}
                         v((t_0,x_0,s_0,y_0)\otimes ({t},{x},s,y))\geq m_0.
$$
                   If $\bar{m}_0=m_0$, then we can take
                   $(\bar{t},\bar{x},\bar{s},\bar{y})=(t_0,x_0,s_0,y_0)$ and finish the procedure.
                   Otherwise there exists $(t_1,x_1,s_1,y_1)\in (t_0,T]\times
                    {\cal{D}}_{\bar{Q}}\times(s_0,T]\times
                    {\cal{D}}_{\bar{Q}}$, such that $(t_1,x_1,s_1,y_1)=(t_0,x_0,s_0,y_0)\otimes
                   (t_1,x_1,s_1,y_1) $ and
$$
                  m_1:=v(t_1,x_1,s_1,y_1)\geq \frac{m_0+\bar{m}_0}{2}.
$$
                  We set
$$
                         \bar{m}_1:=\sup_{(t,x,s,y)\in [t_1,T]\times
                         {\cal{D}}_{\bar{Q}}\times [s_1,T]\times
                         {\cal{D}}_{\bar{Q}}}
                         v((t_1,x_1,s_1,y_1)\otimes ({t},{x},s,y))\geq m_1.
$$
            If $\bar{m}_1=m_1$, then we can take
            $(\bar{t},\bar{x},\bar{s},\bar{y})=(t_1,x_1,s_1,y_1)$ and finish the
            procedure. Otherwise we can find, for $i=2,3,\cdots$,
            $(t_i, x_i,s_i, y_i)\in (t_{i-1},T]\times {\cal{D}}_{\bar{Q}}\times (s_{i-1},T]\times {\cal{D}}_{\bar{Q}}$
            such that $(t_i,x_i,s_i, y_i)=(t_{i-1},x_{i-1},s_{i-1},y_{i-1})\otimes
           (t_i,x_i,s_i, y_i)$, $v(t_i,x_i,s_i, y_i)\geq v(t_i,x_i+e1_0(\cdot),s_i,y_i+l1_0(\cdot))$, for all $e,l$ such
           that $e+x_i(0),l+y_i(0)\in \bar{Q}$ and
$$
                    m_i:=u(t_i, x_i,s_i,y_i)\geq
                    \frac{m_{i-1}+\bar{m}_{i-1}}{2},
$$
$$
                          \bar{m}_i:=\sup_{(t,x,s,y)\in [t_i,T]\times
                         {\cal{D}}_{\bar{Q}}[s_i,T]\times
                         {\cal{D}}_{\bar{Q}}}
                         u((t_i,x_i,s_i, y_i)\otimes ({t},{x},s,y))\geq m_i,
$$
              and continue this procedure till the first time when
              $\bar{m}_i=m_i$ and the finish the proof by setting
              $(\bar{t},\bar{x},\bar{s},\bar{y})=(t_i,x_i,s_i,y_i)$. For the last case in which
              $\bar{m}_i>m_i$ for all $i=1,2,\cdots$, we have $t_i\uparrow \bar{t}\in
              [0,T], s_i\uparrow \bar{s}\in[0,T]$. Then we can find $\bar{x},\bar{y}\in {\cal{D}}_{\bar{Q}}$ such
              that $(\bar{t},\bar{x},\bar{s},\bar{y})=(t_i,x_i,s_i,y_i)\otimes (\bar{t},\bar{x},\bar{s},\bar{y})$.
              We can choose $\bar{x}(0),\bar{y}(0)\in \bar{Q} $ such that  $u(\bar{t},\bar{x},\bar{s},\bar{y})
              \geq u(\bar{t},\bar{x}+x1_0(\cdot),\bar{s},\bar{y}+y1_0(\cdot))$,
              for all $x,y$ such that $x+\bar{x}(0),y+\bar{y}(0)\in \bar{Q}$. Since
$$
                   \bar{m}_{i+1}-m_{i+1}\leq
                   \bar{m}_i-\frac{\bar{m}_i+m_i}{2}=\frac{\bar{m}_i-m_i}{2},
$$
           thus there exists  $\bar{m}\in (m_0,\bar{m}_0)$, such
           that $\bar{m}_i\downarrow \bar{m}$ and ${m}_i\uparrow
           \bar{m}$. By the definitions of $\bar{x}$ and $\bar{y}$,
           we get  $x_i(s)\rightarrow\bar{x}(s)$ and $y_i(s)\rightarrow\bar{y}(s)$ for almost all
           $s\in [-\tau,0]$ and there exist two subsequences of
           $x_i(0)$ and $y_i(0)$ still denoted by themselves such that $x_i(0)\rightarrow \bar{a} \in
           \bar{Q}$ and $y_i(0)\rightarrow \bar{b} \in
           \bar{Q}$, respectively.\\
           Thus, by (4.1) we get that
$$
           \bar{m}=\lim_{i\rightarrow\infty}m_i=\lim_{i\rightarrow\infty}v(t_i,x_i,s_i,y_i)\leq
           v(\bar{t},\bar{x},\bar{s},\bar{y}).
$$
            We can claim that (4.2) holds for this
           $(\bar{t},\bar{x},\bar{s},\bar{y})$. Indeed, otherwise there exist $(t,x,s,y)\in (\bar{t},T]\times
           {\cal{D}}_{\bar{Q}}\times (\bar{s},T]\times
           {\cal{D}}_{\bar{Q}}$ and $\delta>0$  with $(t,x,s,y)=(\bar{t},\bar{x},\bar{s},\bar{y})\otimes
           (t,x,s,y)$, such that
$$
                      v((\bar{t},\bar{x},\bar{s},\bar{y})\otimes
           (t,x,s,y))\geq v(\bar{t},\bar{x},\bar{s},\bar{y})+\delta \geq \bar{m}+\delta,
$$
then the following contradiction is induced:
$$
 v((\bar{t},\bar{x},\bar{s},\bar{y})\otimes
           (t,x,s,y))= v(({t}_i,{x}_i,s_i,y_i)\otimes
           (t,x,s,y))\leq \bar{m}_i\rightarrow \bar{m}.
$$
The proof is completed.  \ \ $\Box$
                                  % consider the viscosity solution of HJB equations (3.13). In fact, we shall get that the value
%                                   function $V$ defined by (3.3) is a viscosity solution of HJB equations (3.13).
                 %                  For simplicity, we define
%%\begin{eqnarray*}
%%\Phi&=&\{\varphi\in C^{1} ([0,T]\times  {\cal{D}})\cap
%%             D({\cal{S}})| \ \exists\ \varphi_0\in C^{1}([0,T]\times R^d\times H),\\
%%                             &&~ \ \mbox{such that }  \varphi(t,x)=\varphi_0(t,x(0),x),  \ \ \ \forall (t,x)\in[0,T]\times {\cal{D}}
%%\};
%%\end{eqnarray*}
%\begin{eqnarray*}
%            {\cal{G}}&=&\{g\in C^{1} ( [0,T]\times H)|\ \exists\ g_0\in C^{0,1}([0,T]\times R)\ \mbox{with} \ g'_0(t,r)\geq0,\ g'_0(t,0)=0,\\
%                             &&~ \ \mbox{such that }  g(t,x)=g_0(t,|x|_H^2),  \ \ \ \forall (t,x)\in[0,T]\times {\cal{D}}\},
%\end{eqnarray*}
%                  where  $g'_0$ denotes first derivative with respect to the second variable of $g_0$.
\par
                                   From the above lemma, we can now give the following definition of viscosity solution:
\begin{definition}
                             $w\in C([0,T]\times {\cal{D}})$ is called a
                             viscosity subsolution (supersolution)
                             of (3.8) if the terminal condition $w(T,x)\leq \phi(x(0))$(resp. $w(T,x)\geq \phi(x(0))$)
                             is satisfied and for every bounded open subset $Q$ of $R^d$ and $\varphi\in \Phi$, whenever the function
                             $w- \varphi $ (resp.  $w+\varphi$) satisfies
$$
                         (w- \varphi)(s,z)=\sup_{(t,x)\in [s,T]\times
                         {\cal{D}}_{\bar{Q}}}
                         (w- \varphi)((s,z)\otimes
                         ({t},{x})),
$$
$$
                       (\mbox{respectively},\ \   (w+\varphi)(s,z)=\inf_{(t,x)\in [s,T]\times
                         {\cal{D}}_{\bar{Q}}}
                         (w+\varphi)((s,z)\otimes
                         ({t},{x}))),
$$
                   where $(s,z)\in [0,T)\times {\cal{D}}_{\bar{Q}}$ and $z(0)\in Q$, we have
\begin{eqnarray*}
                           \varphi_{t}(s,z)+{\cal{S}}(\varphi)(s,z)
                                     +H(s,z,\nabla_{x}\varphi(s,z))\geq0,
\end{eqnarray*}
\begin{eqnarray*}
                          (\mbox{respectively},
                          -\varphi_{t}(s,z)-{\cal{S}}(\varphi)(s,z)+H(s,z,-\nabla_{x}\varphi(s,z))\leq0).
\end{eqnarray*}
                                $w\in C([0,T]\times {\cal{D}})$ is said to be a
                             viscosity solution of (3.8) if it is
                             both a viscosity subsolution and a viscosity
                             supersolution.
\end{definition}
 \begin{remark}
  \begin{description}
  \rm{
  \item{(i)}  A viscosity solution of the HJB equation (3.8)  is a
                       classical solution if it furthermore lies in $\Phi$.
 \item{(ii)}    In the classical uniqueness proof of viscosity solution to HJB equation in infinite
                dimensions, the weak compactness of separable
                Hilbert spaces is used  ( see  \cite{li}). In our case, the HJB
                equation is defined on space $\cal{D}$, which
                doesn't have weak compactness. For the sake of the
uniqueness proof,
   our new notion of viscosity solution is enhanced. At the same time, our modification doesn't lead to additional difficulty in the existence
proof.
\item{(iii)}  Assume that  the coefficient $F(t,x,y,u)=\overline{F}(t,x,u),\  (t,x,y,u)\in [0,T]\times R^d\times R\times U$ and $b=0$.
           Let  function $V(t,x):[0,T]\times R^d\rightarrow R$ be a
             viscosity solution to (3.8) as a  functional $V(t,x):[0,T]\times
             {\cal{D}}\rightarrow R$. Then $V$ is also a classical viscosity
solution as a function of time and state.}
\end{description}
\end{remark}
         We conclude this section with   the existence result on viscosity solution.
\begin{theorem}
                          Suppose that Hypothesis 2.1  and
                          Hypothesis 3.1  hold true. Then the value
                          function $V(t,x)$ defined by (3.3) is a
                          viscosity solution of (3.8).
\end{theorem}
  \par
{\bf  Proof}. \ \        First, for every bounded open subset
                  $Q\subset R^d$ we let $\varphi\in \Phi$
%C^{1}
%([0,T)\times  {\cal{D}})\cap            D({\cal{S}})
 such that $$
                         0=(V- \varphi )(t,x)=\sup_{(s,y)\in [t,T]\times
                         {\cal{D}}_{\bar{Q}}}
                         (V- \varphi)((t,x)\otimes
                         ({s},{y})),
$$
 where $(t,x)\in [0,T)\times {\cal{D}}_{\bar{Q}}$ and $x(0)\in Q$. % Without loss of generality, by adding some
%                             constants to $\varphi$ and $g$ if necessary, we assume that
%$$
%                         g(t,x)=0, \ \ \ \ V(t,x)=\varphi(t,x).
%$$
                   Then, for fixed $u\in U$ and $t\leq s < T$, by the dynamic programming principle (Theorem 3.3), we get
                   that, for $s$ small enough,
\begin{eqnarray*}
                 \varphi(t,x)=V(t,x)
                 \leq \int_{t}^{s}q(\sigma,X^u(\sigma),u)d\sigma+V(s,X^u_s)\leq\int_{t}^{s}q(\sigma,X^u(\sigma),u)d\sigma+\varphi(s,X^u_s).
\end{eqnarray*}
                    Thus,
  \begin{eqnarray*}
                 0  \leq\frac{1}{s-t}\int_{t}^{s}q(\sigma,X^u(\sigma),u)d\sigma
                            +\frac{1}{s-t}{[}\varphi(s,X^u_s)-\varphi(t,x){]}.
\end{eqnarray*}
             Now, applying Lemma 3.5, we show that
\begin{eqnarray*}
                 0  \leq
                     q(t,x(0),u)+ \varphi_{t}(t,x)+{\cal{S}}(\varphi)(t,x)+\langle{\nabla_{x}\varphi(t,x)},
                     [F(t,x(0),(a,x)_{H},u)+b(t)x(-\tau)]1_0(t)\rangle.
\end{eqnarray*}
   Taking the minimum in $u\in U$, we get  that $V$ is a viscosity subsolution of (3.8).
\par
      Next,  for every bounded open subset
                  $Q\subset R^d$, we let     $\varphi\in\Phi$
%([0,T)\times  {\cal{D}})\cap            D({\cal{S}})$
 such that
$$
           0= (V+\varphi)(t,x)=\inf_{(s,y)\in
                [t,T]\times
                         {\cal{D}}_{\bar{Q}}}
                         (V+\varphi)((t,x)\otimes
                         ({s},{y})),
$$
 where $(t,x)\in [0,T)\times {\cal{D}}_{\bar{Q}}$ and $x(0)\in Q$.
 %We also assume that
%$$
%                         g(t,x)=0, \ \ \ \ V(t,x)=-\varphi(t,x).
%$$
   For any $\varepsilon>0$ and $s>t$, by (3.7), one can find a control  ${u}^\varepsilon(\cdot)\equiv u^{\varepsilon,s}(\cdot)\in {\cal{U}}[t,T]$ such
   that, for $s$ small enough,
\begin{eqnarray*}
    \varepsilon(s-t)
    &\geq&\int_{t}^{s}q(\sigma,X^{{u}^\varepsilon}(\sigma),{u}^\varepsilon(\sigma))d\sigma
                                       +V(s,X^{{u}^\varepsilon}_s)-V(t,x)\\
                      &\geq&\int_{t}^{s}q(\sigma,X^{{u}^\varepsilon}(\sigma),{u}^\varepsilon(\sigma))d\sigma
                                -\varphi(s,X^{{u}^\varepsilon}_s)+\varphi(t,x).
\end{eqnarray*}
                 Then, by Lemma 3.5, we obtain that
\begin{eqnarray*}
                           \varepsilon&\geq& \frac{1}{s-t}\int_{t}^{s}q(\sigma,X^{{u}^\varepsilon}(\sigma),{u}^\varepsilon(\sigma))d\sigma
                                -\frac{\varphi(s,X^{{u}^\varepsilon})-\varphi(t,x)}{s-t}\\
           &\geq&-\varphi_t(t,x)-{\cal{S}}(\varphi)(t,x)+\frac{1}{s-t}\int_{t}^{s}q(t,x(0),{u}^\varepsilon(\sigma)) \\
                                &&-\langle{\nabla_{x}\varphi(t,x)}, [F(t,x(0),(a,x)_{H},u^\varepsilon(\sigma))+b(t)x(-\tau)]1_0(t)\rangle d\sigma+o(1)\\
            &\geq&-\varphi_t(t,x)-{\cal{S}}(\varphi)(t,x)+\inf_{u\in U}[q(t,x(0),{u}) \\
            && -\langle{\nabla_{x}\varphi(t,x)}, [F(t,x(0),(a,x)_{H},u)+b(t)x(-\tau)]1_0(t)\rangle]+o(1).
\end{eqnarray*}
Letting $s\downarrow t$ and $\varepsilon\rightarrow0$ we show that
\begin{eqnarray*}
            0\geq-\varphi_t(t,x)-{\cal{S}}(\varphi)(t,x)+\inf_{u\in U}[q(t,x(0),{u})
             -\langle{\nabla_{x}\varphi(t,x)},
             [F(t,x(0),(a,x)_{H},u)+b(t)x(-\tau)]1_0(t)\rangle].
\end{eqnarray*}
    Therefore, $V$ is also a viscosity supsolution of (3.8). This completes the proof of Theorem 4.4.\ \ $\Box$
%%%%%%%%%%%%%%%%%%%%%%%%%%%%%%%%%%%%%%%%%%%%%%%%%%%%%%%%%%%%%%%%%%%%%%%%%%%%%%%%%%%%%%%%%%%%%%%%%%%%%%%%%%%%%%%%%%%%%%%%%

%%%%%%%%%%%%%%%%%%%%%%%%%%%%%%%%%%%%%%%%%%%%%%%%%%%%%%%%%%%%%%%%%%%%%%%%%%%%%%%%%%%%%%%%%%%%%%%%%%%%%%%%%%%%%%%%%%%%%%%%%

\section{Viscosity solution of HJB equations: Uniqueness theorem.}
\par
                  This section is devoted to a  proof of uniqueness of the viscosity
                   solution to (3.8). This result, together with
                   those in the previous section, will give a
                   characterization for the value function of
                   optimal control problem (3.1) and (3.2).
\par
            We are now state the main result of this section.
\begin{theorem}  Suppose that Hypothesis 2.1 and
                          Hypothesis 3.1 hold true.
                        Let $W$ (resp. $V$) be a viscosity subsolution  (resp. supsolution) of (3.8) and  there exists a constant $\Lambda>0$
                        such that, for $(t,x), (s,y)\in [0,T]\times {\cal{D}}$,
   $$
                                    |W(t,x)|\vee|V(t,x)|\leq \Lambda (1+|x(0)|+|x|_H),   \eqno(5.1)
   $$
                         and %assume if $x,y$ are bounded in ${\cal{D}}$, we have
%                         that, for some constant $C>0$,
\begin{eqnarray*}
                        &&|W(t,x)-W(s,y)|\vee|V(t,x)-V(s,y)|
                  \leq
                        \Lambda(1+|x(0)|+|x|_H+|y(0)|+|y|_H)\\
                        &&~~~~~~~~~~~~~~~~~~~~~~~~~~~~~~~~~~~\times\bigg{(}|x(0)-y(0)|+|s-t|^{\frac{1}{2}}
                  +\sup_{l\in[-\tau,0]}\bigg{|}\int^{0}_{l}x(\theta)-y(\theta)d\theta\bigg{|}\bigg{)}.
                  \ \ \
                  (5.2)
\end{eqnarray*}
 %$$
%                            |W(t,x)-W(s,y)|+|V(t,x)-V(s,y)|\rightarrow0\
%                            \mbox{as}\
%                            |s-t|+|x(0)-y(0)|+|x-y|\rightarrow0.\eqno(5.2)
% $$
                         Then $W\leq V$.
\end{theorem}
\par
          From this theorem, the viscosity solution to HJB equation (3.8) can characterizes the   value function $V(t,x)$ of our optimal control
          problem (3.1) and (3.2)  as following:
\begin{theorem}  Let Hypothesis 2.1 and
                          Hypothesis 3.1 hold true. Then  the value function $V$ defined by (3.3) is the
                            unique viscosity solution of (3.8).
\end{theorem}
\par
{\bf  Proof}. \ \
               By Theorem 4.4, we know that $V$ is a viscosity solution of (3.8). Thus, our conclusion follows from Theorem 3.2 and Theorem 5.1.\ \ $\Box$
   \par We are now in a position of showing the proof the Theorem
   5.1.
 We first note that for $\delta>0$, the function
                    defined by $\tilde{W}:=W-\frac{\delta}{t}$ is a subsolution
                    of
 \begin{equation*}
\begin{cases}
\frac{\partial}{\partial t} \tilde{W}(t,x)+{\mathcal
                        {S}}(\tilde{W})(t,x)+ H(t,x,\nabla_x \tilde{W}(t,x))= \frac{\delta}{t^2}, \ \  t\in
                               [0,T],\ \ x\in {\mathcal{D}}, \\
 \tilde{W}(T,x)=\phi(x(0)).
\end{cases}
\end{equation*}
%$$
%          \cases{\frac{\partial}{\partial t} \tilde{W}(t,x)+{\mathcal
%                        {S}}(\tilde{W})(t,x)+ H(t,x,\nabla_x \tilde{W}(t,x))= \frac{\delta}{t^2}, \ \  t\in
%                               [0,T],\ \ x\in {\mathcal{D}}, \cr
%         \tilde{W}(T,x)=\phi(x(0)).\cr}
%$$
                Since $W\leq V$ follows from $\tilde{W}\leq V$ in
                the limit $\delta\downarrow0$, it suffices to prove
                the theorem under the additional assumption:
$$
                 \frac{\partial}{\partial t} {W}(t,x)+{\mathcal
                        {S}}({W})(t,x)+ H(t,x,\nabla_x {W}(t,x))\geq c,\ \ c:=\frac{\delta}{T^2}, \ \  t\in
                               [0,T],\ \ x\in {\mathcal{D}}.
$$
\par
   {\bf  Proof  of Theorem 5.1}. \ \   The proof of this theorem  is rather long. Thus, we split it into several
        steps.
\par
            $Step\  1.$ Definition of auxiliary functions and sets.
\par
We only need to prove that $W(t,x)\leq V(t,x)$ for all $(t,x)\in
[T-\bar{a},T)\times
        {\cal{D}}$.
        Here
        $$\bar{a}=\frac{1}{8(1+L)^2{\bar{C}}^2}\wedge\frac{\tau}{2},\
        \ \
\bar{C}={{1+\tau|a|_{W^{1,2}}+\sup_{s\in[0,T]}|b(s)|}}.$$
         Then repeat the same procedure for cases
        $[T-i\bar{a},T-(i-1)\bar{a})$.  To this end we assume to the
        contrary that there exists $(\tilde{t},\breve{x})\in [T-\bar{a},T)\times
        {\cal{D}}$, such that
        $2\tilde{m}:=W(\tilde{t},\breve{x})-V(\tilde{t},\breve{x})>0$.
        Since Lipschitz continuous functions  are dense in $H$, by (5.2) there exist a
        Lipschitz continuous function $\tilde{y}$ and $\tilde{a}\in R^d$ such that
        $W(\tilde{t},\tilde{x})-V(\tilde{t},\tilde{x})>\tilde{m}$,
        where $\tilde{x}=\tilde{y}+\tilde{a}1_0(\cdot)$.
\par
         First, let $\varepsilon >0$ be a small number such that
 $$
 W(\tilde{t},\tilde{x})-V(\tilde{t},\tilde{x})-2\varepsilon \frac{\mu T-\tilde{t}}{\mu
 T}(|\tilde{x}|_H^2+|\tilde{x}(0)|^2)>\frac{\tilde{m}}{2},
 $$
      and
 $$
                          \frac{9L^2\varepsilon}{8\mu T(1+L)^2}\leq\frac{c}{2}, \eqno(5.3)
$$
             where
$$
            \mu=1+\frac{1}{4T(1+L)^2\bar{C}^2}.
$$
                      Next, for every $\alpha>0$ we define for any $(t,x,s,y)\in [0,T]\times {\cal{D}}\times [0,T]\times  {\cal{D}}$,
\begin{eqnarray*}
                 \Psi(t,x,s,y)&=&W(t,x)-V(s,y)-\frac{\alpha}{2}d(t,s,x,y)\\
                 &&-\varepsilon\frac{\mu T-t}{\mu T}(|x|_H^2+|x(0)|^2)
            -\varepsilon\frac{\mu T-s}{\mu T}(|y|_H^2+|y(0)|^2),
\end{eqnarray*}
 where
 $$
            d(t,x,s,y)=|x(0)-y(0)|^2+|x-y|_B^2+|s-t|^2.
 $$
        %Let $M$ be a large number such that $
%        \tilde{x}\in {\cal{D}}_{{Q}^M}$ and for all
%        $|x(0)|=|y(0)|=M$, $x,y\in {\cal{D}}_{\overline{Q^M}} $,
%  We define
% $$
%                            \Psi(t,x,s,y)<  M_\alpha-1:= \sup_{t,s\geq \tilde{t}; x,y\in
%                              {\cal{D}}_{\overline{Q^M}}} \Psi((\tilde{t},\tilde{x},\tilde{t},\tilde{x})\otimes (t,x,s,y))-1,
% $$
  Finally, for every $M>0$ satisfying $\tilde{x}\in {\cal{D}}_{\overline{Q^M}}$, we define
 $$
                              M_\alpha:= \sup_{t,s\geq \tilde{t}; x,y\in
                              {\cal{D}}_{\overline{Q^M}}} \Psi((\tilde{t},\tilde{x},\tilde{t},\tilde{x})\otimes (t,x,s,y)),
 $$
        and
 $$
                             M_\alpha \geq M_*:=\sup_{t\geq \tilde{t}; x\in
                              {\cal{D}}_{\overline{{{Q}^M}}}} \Psi((\tilde{t},\tilde{x},\tilde{t},\tilde{x})\otimes (t,x,t,x))\geq
                              \frac{\tilde{{m}}}{2},
 $$
         where
$$
                   Q^M:=\{(x_1,x_2,\ldots,x_d)|\
                   |x_1|^2+|x_2|^2+\cdots+|x_d|^2<M^2\}.
$$
 %Let $\alpha$ be a large number such that
% $$
%                           \frac{1}{\alpha^2}+\rho(\frac{2}{\alpha}(\frac{1}{\alpha^2}+C-M_*))\leq
%                           \frac{1}{4}(\tilde{m}\wedge c),
% $$
                 %where $C:=2\Lambda(1+M+\tau^\frac{1}{2}M)$.
\par
$Step\ 2.$ Properties of $\Psi(t,x,s,y)$.
\par
      By the definition of $M_\alpha$, we can fix $(\bar{t},\bar{x}), (\bar{s},\bar{y})\in [\tilde{t},T]\times {\cal{D}}_{{\overline{Q^M}}}$ satisfying
        $$(\bar{t},\bar{x},\bar{s},\bar{y})=(\tilde{t},\tilde{x},\tilde{t},\tilde{x})\otimes
        (\bar{t},\bar{x},\bar{s},\bar{y}),\ \
        \Psi(\tilde{t},\tilde{x},\tilde{t},\tilde{x})\leq \Psi(\bar{t},\bar{x},\bar{s},\bar{y})
        \ \ \mbox{and}\ \ \Psi(\bar{t},\bar{x},\bar{s},\bar{y})+\frac{1}{\alpha}>M_\alpha.$$
Now we can apply Lemma 4.1 to find $(\hat{{t}},\hat{{x}}),
(\hat{{s}},\hat{{y}})\in [T-\bar{a},T]\times
            {\cal{D}}_{\overline{{Q}^M}}$   satisfying $(\hat{{t}},\hat{{x}},\hat{{s}},\hat{{y}})
            =(\bar{t},\bar{x},\bar{{s}},\bar{{y}})\otimes
            (\hat{{t}},\hat{{x}},\hat{{s}},\hat{{y}})$ %and $(\hat{{s}},\hat{{y}})=(\bar{s},\bar{y})\otimes
%            (\hat{{s}},\hat{{y}})$
          with $\Psi(\hat{t},\hat{x},\hat{s},\hat{y})\geq \Psi(\bar{t},\bar{x},\bar{s},\bar{y})\geq \Psi(\tilde{t},\tilde{x},\tilde{t},\tilde{x})$ such that
 $$
                   \Psi(\hat{t},\hat{x},\hat{s},\hat{y})\geq  \Psi((\hat{{t}},\hat{{x}},\hat{{s}},\hat{{y}})\otimes (t,x,s,y)),\  t\geq \hat{t}, s\geq
                              \hat{s}, x,y \in
                              {\cal{D}}_{\overline{{Q}^M}}.
 $$
             We should note that the point
             $(\hat{t},\hat{x},\hat{s},\hat{y})$ depends on $\bar{t},\bar{x},\bar{s},\bar{y},\alpha,
             M$.
 \par
        $Step \ 3.$ For fixed $
             M$, there exists a  subsequence of $\alpha$ still
        denoted
               by itself such that
$$
                         \frac{\alpha}{2}d(\hat{t},\hat{x},\hat{s},\hat{y})
                         \leq\frac{1}{\alpha}+|W(\hat{t},\hat{x})-W(\hat{s},\hat{y})|
                                   +|V(\hat{t},\hat{x})-V(\hat{s},\hat{y})|\rightarrow0 \ \mbox{as}\ \alpha\rightarrow+\infty,\eqno(5.4)
$$
and
$$
                            \alpha|b(\hat{t})\hat{x}(-\tau)-b(\hat{s})\hat{y}(-\tau)|^2
                            \rightarrow0\
                            \mbox{as}\
                            \alpha\rightarrow+\infty.\eqno(5.5)
$$
      Let us show the above.  We can check that
 \begin{eqnarray*}
                         ~~~~~~~~~~~~~~~~~~&&\frac{\alpha}{2}d(\hat{t},\hat{x},\hat{s},\hat{y})+\varepsilon\frac{\mu T-\bar{t}}{\mu T}(|\hat{x}|_H^2+|\hat{x}(0)|^2)
                             +\varepsilon\frac{\mu T-\bar{s}}{\mu T}(|\hat{y}|_H^2+|\hat{y}(0)|^2)\\
                         &\leq& \frac{1}{\alpha}+W(\hat{t},\hat{x})-V(\hat{s},\hat{y})-M_\alpha\leq\frac{1}{\alpha}+W(\hat{t},\hat{x})-V(\hat{s},\hat{y})-M_*\\
                         &\leq&\frac{1}{\alpha}+C-M_*,\ \ \ \ \ ~~~~~~~~~~~~~~~~~~~~~~~~~~~~~~~~~~~~~~~~~~~~~~~~~~~~~~~~~~~~~~~
                                               ~~~~~~~~~ \ (5.6)
 \end{eqnarray*}
 where $C:=2\Lambda(1+M+\tau^\frac{1}{2}M)$.
 We also have that
 \begin{eqnarray*}
                         2M_*&\leq&\frac{2}{\alpha}+W(\hat{t},\hat{x})-W(\hat{s},\hat{y})+W(\hat{s},\hat{y})-V(\hat{s},\hat{y})
                                    +W(\hat{t},\hat{x})-V(\hat{t},\hat{x})+V(\hat{t},\hat{x})-V(\hat{s},\hat{y})\\
                              &&-{\alpha}d(\hat{t},\hat{x},\hat{s},\hat{y})
                              -2\varepsilon\frac{\mu T-\bar{t}}{\mu T}(|\hat{x}|_H^2+|\hat{x}(0)|^2)
                             -2\varepsilon\frac{\mu T-\bar{s}}{\mu T}(|\hat{y}|_H^2+|\hat{y}(0)|^2)\\
                         &\leq&\frac{2}{\alpha}+|W(\hat{t},\hat{x})-W(\hat{s},\hat{y})|
                                   +|V(\hat{t},\hat{x})-V(\hat{s},\hat{y})|
                                +2M_*-{\alpha}d(\hat{t},\hat{x},\hat{s},\hat{y}).
 \end{eqnarray*}
  Thus
 $$
                         \frac{\alpha}{2}d(\hat{t},\hat{x},\hat{s},\hat{y})
                         \leq\frac{1}{\alpha}+|W(\hat{t},\hat{x})-W(\hat{s},\hat{y})|
                                   +|V(\hat{t},\hat{x})-V(\hat{s},\hat{y})|.\eqno(5.7)
$$
                 By the definition of $d$ and (5.6), we get that
 $
                          |\hat{x}-\hat{y}|_B\rightarrow0\
                          \mbox{as}\ \alpha\rightarrow+\infty$. We
                          note that $|\hat{x}|_{{\cal{D}}}\vee|\hat{y}|_{{\cal{D}}}\leq
                          M$.
        Then, from the definition of $B$, it follows that there exists a subsequence of $\alpha$ still
        denoted
               by itself, such that
                $$\sup_{l\in[-\tau,0]}\bigg{|}\int^{0}_{l}\hat{x}(\theta)-\hat{y}(\theta)d\theta\bigg{|}\rightarrow0\
                          \mbox{as} \ \alpha\rightarrow+\infty.
                          $$
                           Combining  (5.2) and (5.7) we see
                           that (5.4) holds.
 %$$
%                         \frac{\alpha}{2}d(\hat{t},\hat{x},\hat{s},\hat{y})
%                         \leq\frac{1}{\alpha}+|W(\hat{t},\hat{x})-W(\hat{s},\hat{y})|
%                                   +|V(\hat{t},\hat{x})-V(\hat{s},\hat{y})|\rightarrow0 \ \mbox{as}\ \alpha\rightarrow+\infty.
%$$
  On the other hand, by  $\tilde{x}=\tilde{y}+\tilde{a}1_0(\cdot)$,  $0<\bar{a}\leq
     \frac{\tau}{2}$, $(\hat{{t}},\hat{{x}},\hat{{s}},\hat{{y}})=(\bar{t},\bar{x},\bar{s},\bar{y})\otimes
            (\hat{{t}},\hat{{x}},\hat{{s}},\hat{{y}})$ and $(\bar{t},\bar{x},\bar{s},\bar{y})
            =(\tilde{t},\tilde{x},\tilde{t},\tilde{x})\otimes (\bar{t},\bar{x},\bar{s},\bar{y})$, we have that
 $$
 \hat{x}(-\tau)=\tilde{y}(\hat{t}-\tau-\tilde{t}),\ \ \
 \hat{y}(-\tau)=\tilde{y}(\hat{s}-\tau-\tilde{t}).
 $$
               Since $b$ and $\tilde{y}$ are Lipschitz continuous,  there
                exists a constant $N>0$ such that
$$
                            \alpha|b(\hat{t})\hat{x}(-\tau)-b(\hat{s})\hat{y}(-\tau)|^2
                            =\alpha|b(\hat{t})\tilde{y}(\hat{t}-\tau-\tilde{t})-b(\hat{s})\tilde{y}(\hat{s}-\tau-\tilde{t})|^2\leq N\alpha |\hat{t}-\hat{s}|^2.
                            %\rightarrow0\
%                            \mbox{as}\
%                            \alpha\rightarrow+\infty.
$$
              Then,  (5.5) follows from (5.6).
\par
                   $Step\ 4.$  There exist $N,M>0$ such that
                   $\hat{t},\hat{s}\in [\tilde{t},T)$ and  $\hat{x}(0),\hat{y}(0)\in
                 Q^M$ for all $\alpha\geq N$.
\par
           First, we  note that, for any $\alpha>0$, there exists a $M>0$ be large enough such that
\begin{eqnarray*}
                     &&\Psi(\tilde{t},\tilde{x},\tilde{t},\tilde{x})=W(\tilde{t},\tilde{x})-V(\tilde{t},\tilde{x})
                 -2\varepsilon\frac{\mu T-\tilde{t}}{\mu T}(|\tilde{x}|_H^2+|\tilde{x}(0)|^2)\\
            &>&W(t,x)-V(s,y)-\varepsilon\frac{\mu T-t}{\mu T}M^2
            -\varepsilon\frac{\mu T-s}{\mu T}M^2\\
            &\geq&\Psi({t},{x},{s},{y}),
\end{eqnarray*}
          where $(t,x,s,y)\in [0,T]\times {\cal{D}}\times [0,T]\times
          {\cal{D}}$ and $|x(0)|=|y(0)|=M$.
            Therefore, for this $M>0$, we have that  $\hat{x}(0),\hat{y}(0)\in
                 Q^M$ for every $\alpha>0$.
\par
  Next, by (5.4), we can let $N>0$ be a large number such that
$$
                         \frac{\alpha}{2}d(\hat{t},\hat{x},\hat{s},\hat{y})\leq\frac{1}{\alpha}+|W(\hat{t},\hat{x})-W(\hat{s},\hat{y})|
                                   +|V(\hat{t},\hat{x})-V(\hat{s},\hat{y})|
                         \leq
                         \frac{1}{4}(\tilde{m}\wedge c),
$$
               for all $\alpha\geq N$.
            Moreover, we have $\hat{t},\hat{s}\in [\tilde{t},T)$ for all $\alpha\geq N$. Indeed, if say $\hat{s}=T$, then we will deduce the following contradiction:
 \begin{eqnarray*}
                         &&\frac{\tilde{m}}{2}\leq M_*\leq M_\alpha \leq\frac{1}{\alpha}+\Psi(\hat{t},\hat{s},\hat{x},\hat{y})\\
                         &\leq&\frac{1}{\alpha}+ W(\hat{t},\hat{x})-W(\hat{s},\hat{y})+W(\hat{s},\hat{y})-V(\hat{s},\hat{y})\\
                         &\leq& \frac{1}{\alpha}+|W(\hat{t},\hat{x})-W(\hat{s},\hat{y})|\leq \frac{\tilde{m}}{4}.
 \end{eqnarray*}
            %Now we can apply Lemma 4.4 to find $(\hat{{t}},\hat{{x}}), (\hat{{s}},\hat{{y}})\in [T-\bar{a},T)\times
%            {\cal{D}}_{\overline{{Q}^M}}$   satisfying $(\hat{{t}},\hat{{x}},\hat{{s}},\hat{{y}})
%            =(\bar{t},\bar{x},\bar{{s}},\bar{{y}})\otimes
%            (\hat{{t}},\hat{{x}},\hat{{s}},\hat{{y}})$ %and $(\hat{{s}},\hat{{y}})=(\bar{s},\bar{y})\otimes
%%            (\hat{{s}},\hat{{y}})$
%          with $\Psi(\hat{t},\hat{x},\hat{s},\hat{y})\geq \Psi(\bar{t},\bar{x},\bar{s},\bar{y})\geq \Psi(\tilde{t},\tilde{x},\tilde{t},\tilde{x})$ such that
% $$
%                   \Psi(\hat{t},\hat{x},\hat{s},\hat{y})\geq  \Psi((\hat{{t}},\hat{{x}},\hat{{s}},\hat{{y}})\otimes (t,x,s,y)),\  t\geq \hat{t}, s\geq
%                              \hat{s}, x,y \in
%                              {\cal{D}}_{\overline{{Q}^M}}.
% $$
\par
          $Step\ 5.$    Completion of the proof.
\par
          From above all,  for the fixed  $N,M>0$ in step 4, we  find
$(\hat{{t}},\hat{{x}}), (\hat{{s}},\hat{{y}})\in [\tilde{t},T)\times
            {\cal{D}}_{\overline{{Q}^M}}$   satisfying $\hat{x}(0),\hat{y}(0)\in
                 Q^M$ for all $\alpha\geq N$ and $(\hat{{t}},\hat{{x}},\hat{{s}},\hat{{y}})
            =(\bar{t},\bar{x},\bar{{s}},\bar{{y}})\otimes
            (\hat{{t}},\hat{{x}},\hat{{s}},\hat{{y}})$ %and $(\hat{{s}},\hat{{y}})=(\bar{s},\bar{y})\otimes
%            (\hat{{s}},\hat{{y}})$
          with $\Psi(\hat{t},\hat{x},\hat{s},\hat{y})\geq \Psi(\bar{t},\bar{x},\bar{s},\bar{y})\geq \Psi(\tilde{t},\tilde{x},\tilde{t},\tilde{x})$ such that
 $$
                   \Psi(\hat{t},\hat{x},\hat{s},\hat{y})\geq  \Psi((\hat{{t}},\hat{{x}},\hat{{s}},\hat{{y}})\otimes (t,x,s,y)),\  t\geq \hat{t}, s\geq
                              \hat{s}, x,y \in
                              {\cal{D}}_{\overline{{Q}^M}}.\eqno(5.8)
 $$
  Then
 $$
 \Psi(\hat{t},\hat{x},\hat{s},\hat{y})\geq  \Psi((\hat{{t}},\hat{{x}})\otimes (t,x),\hat{s},\hat{y}),\  t\geq \hat{t},  x \in
                              {\cal{D}}_{\overline{{Q}^M}}.
 $$
      Thus, by the definition of the viscosity subsolution, we get
         that
\begin{eqnarray*}
                      &&\alpha(\hat{t}-\hat{s})-\frac{\varepsilon}{\mu T}(|\hat{x}|_H^2+|\hat{x}(0)|^2)
                      +\frac{\alpha}{2}( B(\hat{x}-\hat{y}), \hat{x}(0)1_{[-\tau,0]}-\hat{x})_H+\varepsilon\frac{\mu T-\hat{t}}{\mu T}(|\hat{x}(0)|^2-|\hat{x}(-\tau)|^2)\\
                                     &&+H(\hat{t},\hat{x},[{\alpha}(\hat{x}(0)-\hat{y}(0))+2\varepsilon\frac{\mu T-\hat{t}}{\mu T}\hat{x}(0)]1_0(t))
                                     \geq c, \ \ \ \ \ \ \ \ \ \ \ \
                                     \ \ \ \ \ \ \ \ \ \ \ \ \ \ \ \
                                     \ \ \ \ \ \ \ \ \ \ \ \ \ \ \ \
                                     \ \
                                     \ \ (5.9)
\end{eqnarray*}
Also, by (5.8), we have
$$
                   \Psi(\hat{t},\hat{x},\hat{s},\hat{y})\geq  \Psi((\hat{{t}},\hat{{x}},(\hat{{s}},\hat{{y}})\otimes (s,y)),\  s\geq
                              \hat{s}, y \in
                              {\cal{D}}_{\overline{{Q}^M}}.
 $$
 Thus, we obtain
\begin{eqnarray*}
                      &&\alpha(\hat{t}-\hat{s})+\frac{\varepsilon}{\mu T}(|\hat{y}|_H^2+|\hat{y}(0)|^2)
                      -\frac{\alpha}{2}( B(\hat{y}-\hat{x}), \hat{y}(0)1_{[-\tau,0]}-\hat{y})_H-\varepsilon\frac{\mu T-\hat{s}}{\mu T}(|\hat{y}(0)|^2-|\hat{y}(-\tau)|^2)\\
                                    && +H(\hat{s},\hat{y},[-{\alpha}(\hat{y}(0)-\hat{x}(0))-2\varepsilon\frac{\mu T-\hat{s}}{\mu T}\hat{y}(0)]1_0(t))
                                     \leq0. \ \ \ \ \ \ \ \ \ \ \ \
                                     \ \ \ \ \ \ \ \ \ \ \ \ \ \ \ \
                                     \ \ \ \ \ \ \ \ \ \ \ \ \ \ \ \
                                       (5.10)
\end{eqnarray*}
                  Combining (5.9) and (5.10), we obtain
\begin{eqnarray*}
                      ~~~~~&&c+\frac{\varepsilon}{\mu T}(|\hat{x}|_H^2+|\hat{x}(0)|^2+|\hat{y}|_H^2+|\hat{y}(0)|^2) +
                      \varepsilon\frac{\mu T-\hat{t}}{\mu T}|\hat{x}(-\tau)|^2+
                      \varepsilon\frac{\mu T-\hat{s}}{\mu T}|\hat{y}(-\tau)|^2\\
                      &\leq&\frac{\alpha}{2}\langle B(\hat{x}-\hat{y}, \hat{x}(0)1_{[-\tau,0]}-\hat{x}-\hat{y}(0)1_{[-\tau,0]}+\hat{y}\rangle
                       +\varepsilon\frac{\mu T-\hat{t}}{\mu T}|\hat{x}(0)|^2+\varepsilon\frac{\mu T-\hat{s}}{\mu T}|\hat{y}(0)|^2\\
                                     &&+H(\hat{t},\hat{x},[{\alpha}(\hat{x}(0)-\hat{y}(0))+2\varepsilon\frac{\mu T-\hat{t}}{\mu
                                     T}\hat{x}(0)]1_0(t))\\
                                     &&
                                     -H(\hat{s},\hat{y},[{\alpha}(\hat{x}(0)-\hat{y}(0))-2\varepsilon\frac{\mu T-\hat{s}}{\mu
                                     T}\hat{y}(0)]1_0(t)). \ \ \ \ \
                                     \ \ \ \ \ \ \ \ \ \ \ \ \ \ \ \
                                     \ \ \ \ \ \ \ \ \ \ \ \ \ \ \ \
                                     \ \ \ \ \ \ \ \ \ (5.11)
\end{eqnarray*}
                   On the other hand, by simple calculation we obtain
\begin{eqnarray*}
                &&\bigg{|}H(\hat{t},\hat{x},[{\alpha}(\hat{x}(0)-\hat{y}(0))+2\varepsilon\frac{\mu T-\hat{t}}{\mu T}\hat{x}(0)]1_0(t))
                                     -H(\hat{s},\hat{y},[{\alpha}(\hat{x}(0)-\hat{y}(0))-2\varepsilon\frac{\mu T-\hat{s}}{\mu T}\hat{y}(0)]1_0(t))\bigg{|}\\
                &&\leq|\sup_{u\in U}(J_1+J_2)|, \ \ \ \ \
                                     \ \ \ \ \ \ \ \ \ \ \ \ \ \ \ \
                                     \ \ \ \ \ \ \ \ \ \ \ \ \ \ \ \
                                     \ \ \ \ \ \ \ \ \ \ \ \ \ \
                                     \ \ \ \ \ \ \ \ \ \ \ \ \ \ \ \
                                     \ \ \ \ \ \ \ \ \ \ \ \ \ \ \ \
                                     \ \ \ \ \ \ \ \ \ \ (5.12)
\end{eqnarray*}
            where
\begin{eqnarray*}
                                 &&J_1=\langle F(\hat{t},\hat{x}(0),(a,\hat{x})_{H},u)+b(\hat{t})\hat{x}(-\tau),{\alpha}(\hat{x}(0)-\hat{y}(0))
                                              +2\varepsilon\frac{\mu T-\hat{t}}{\mu T}\hat{x}(0)\rangle\\
                                 &&~~~~~~~~-\langle F(\hat{s},\hat{y}(0),(a,\hat{y})_{H},u)+b(\hat{s})\hat{y}(-\tau),{\alpha}(\hat{x}(0)-\hat{y}(0))
                                               -2\varepsilon\frac{\mu T-\hat{s}}{\mu T}\hat{y}(0)\rangle\\
                                 &\leq&\alpha{L}(|\hat{x}(0)-\hat{y}(0)|^2+
                                         |\hat{x}(0)-\hat{y}(0)|[|a|_{W^{1,2}}|\hat{x}-\hat{y}|_B+|b(\hat{t})\hat{x}(-\tau)-b(\hat{s})\hat{y}(-\tau)|+|\hat{t}-\hat{s}|])\\
                                         &&+2\varepsilon L(\frac{\mu T-\hat{t}}{\mu T})|\hat{x}(0)|(1+\sup_{s\in[0,T]}|b(s)||\hat{x}(-\tau)|
                                             +|\hat{x}(0)|+|a|_{W^{1,2}}|\hat{x}|_B)\\
                                         &&+2\varepsilon L(\frac{\mu T-\hat{s}}{\mu T})|\hat{y}(0)|(1+\sup_{s\in[0,T]}|b(s)||\hat{y}(-\tau)|
                                             +|\hat{y}(0)|+|a|_{W^{1,2}}|\hat{y}|_B); \ \ \ \ \ \ \ \ \ \
                                             \ \ \ \
                                     \ \ \ \ \ \ \ \ \ \ \  (5.13)
\end{eqnarray*}
and
$$
                                 J_2=q(\hat{t},\hat{x}(0),u)-
                                 q(\hat{s},\hat{y}(0),u)\leq
                                 L|\hat{x}(0)-\hat{y}(0)|+\rho(|\hat{t}-\hat{s}|,M).\eqno(5.14)
$$
                Combining (5.11)-(5.14), we get
\begin{eqnarray*}
                      c
                      &\leq&-\frac{\varepsilon}{\mu
                      T}(|\hat{x}|_H^2+|\hat{x}(0)|^2+|\hat{y}|_H^2+|\hat{y}(0)|^2)-
                      \varepsilon\frac{\mu T-\hat{t}}{\mu
                      T}|\hat{x}(-\tau)|^2-
                      \varepsilon\frac{\mu T-\hat{s}}{\mu
                      T}|\hat{y}(-\tau)|^2\\
                      &&+\frac{\alpha}{4}( |(\hat{x}-\hat{y})|_B^2+\tau|\hat{x}(0)-\hat{y}(0)|^2)+\varepsilon\frac{\mu T-\hat{t}}{\mu T}|\hat{x}(0)|^2
                                    +\varepsilon\frac{\mu T-\hat{s}}{\mu T}|\hat{y}(0)|^2\\
                                     &&+\alpha{L}(3|\hat{x}(0)-\hat{y}(0)|^2+
                                         |a|^2_{W^{1,2}}|\hat{x}-\hat{y}|^2_B+|b(\hat{t})\hat{x}(-\tau)-b(\hat{s})\hat{y}(-\tau)|^2+|\hat{t}-\hat{s}|^2)\\
                                         &&+2\varepsilon L(\frac{\mu T-\hat{t}}{\mu T})|\hat{x}(0)|(1+\sup_{s\in[0,T]}|b(s)||\hat{x}(-\tau)|
                                             +|\hat{x}(0)|+|a|_{W^{1,2}}|\hat{x}|_B)\\
                                         &&+2\varepsilon L(\frac{\mu T-\hat{s}}{\mu T})|\hat{y}(0)|(1+\sup_{s\in[0,T]}|b(s)||\hat{y}(-\tau)|
                                             +|\hat{y}(0)|+|a|_{W^{1,2}}|\hat{y}|_B)\\
                                         &&
                                         +L|\hat{x}(0)-\hat{y}(0)|+\rho(|\hat{t}-\hat{s}|,M).
\end{eqnarray*}
                         Recalling $\bar{a}=\frac{1}{8(1+L)^2\bar{C}^2}\wedge\frac{\tau}{2}$, $\mu=1+\frac{1}{4T(1+L)^2\bar{C}^2}$
                         and $\hat{t},\hat{s}\in [T-\bar{a},T)$, we show that
\begin{eqnarray*}
                      c
                      &\leq&-\frac{\varepsilon}{4\mu T(1+L)^2}\bigg{(}(\frac{1}{{\bar{C}}}{|\hat{x}(-\tau)|}-\frac{3L}{2}\hat{x}(0))^2
                                         +(2(1+L)|\hat{x}|-\frac{3L}{4(1+L)}|\hat{x}(0)|)^2+(\hat{x}(0)-\frac{3L}{2})^2\bigg{)}\\
                      &&-\frac{\varepsilon}{4\mu T(1+L)^2}\bigg{(}(\frac{1}{{\bar{C}}}|\hat{y}(-\tau)|-\frac{3L}{2}\hat{y}(0))^2
                                         +(2(1+L)|\hat{y}|-\frac{3L}{4(1+L)}|\hat{y}(0)|)^2+(\hat{y}(0)-\frac{3L}{2})^2\bigg{)}\\
                                     &&+\frac{9L^2\varepsilon}{8\mu T(1+L)^2}+\alpha(3L+\tau)(|\hat{x}(0)-\hat{y}(0)|^2
                                       +|b(\hat{t})\hat{x}(-\tau)-b(\hat{s})\hat{y}(-\tau)|^2+|\hat{t}-\hat{s}|^2)\\
                                     &&+\alpha(1+L\bar{C}^2)|\hat{x}-\hat{y}|^2_B+L|\hat{x}(0)-\hat{y}(0)|+\rho(|\hat{t}-\hat{s}|,M)\\
                       &\leq&\frac{9L^2\varepsilon}{8\mu T(1+L)^2}+\alpha(3L+\tau)(|\hat{x}(0)-\hat{y}(0)|^2
                                       +|b(\hat{t})\hat{x}(-\tau)-b(\hat{s})\hat{y}(-\tau)|^2+|\hat{t}-\hat{s}|^2)\\
                                     &&+\alpha(1+L\bar{C}^2)|\hat{x}-\hat{y}|^2_B+L|\hat{x}(0)-\hat{y}(0)|+\rho(|\hat{t}-\hat{s}|,M).
\end{eqnarray*}
                      Letting $\alpha\rightarrow+\infty$, it follows from (5.3) that
\begin{eqnarray*}
                      c
                      \leq\frac{c}{2},
\end{eqnarray*}
which induces a contradiction. The proof is completed. \ \ $\Box$

\end{document}